\documentclass[11pt,reqno]{amsart}
\allowdisplaybreaks[4]
\usepackage{amssymb}
\usepackage{anysize}
\usepackage{hyperref}
\usepackage{extarrows}

\numberwithin{equation}{section}
\newtheorem{thm}{Theorem}[section]
\newtheorem{lem}[thm]{Lemma}

\newtheorem{rem}[thm]{Remark}

\newcommand{\be}{\begin{equation}}
\newcommand{\ee}{\end{equation}}
\newcommand{\bea}{\begin{eqnarray*}}
\newcommand{\eea}{\end{eqnarray*}}

\makeatletter

\newcommand{\Rmnum}[1]{\expandafter\@slowromancap\romannumeral #1@}
\makeatother
\makeatletter 
\@addtoreset{equation}{section}
\makeatother  

\marginsize{35mm}{35mm}{38mm}{40mm}

\begin{document}

\title[]{Central limit theorem for partial linear eigenvalue statistics of Wigner matrices}
\author{Zhigang Bao}
\author{Guangming Pan}
\author{Wang Zhou}
\
\thanks{ Z.G. Bao was partially supported by NSFC grant
11071213, NSFC grant 11101362, ZJNSF grant R6090034 and SRFDP grant 20100101110001;
G.M. Pan was partially supported by the Ministry of Education, Singapore, under grant \# ARC 14/11;
   W. Zhou was partially supported by the Ministry of Education, Singapore, under grant \# ARC 14/11,  and by a grant
R-155-000-116-112 at the National University of Singapore.}

\address{Department of Mathematics, Zhejiang University, P. R. China}
\email{zhigangbao@zju.edu.cn}

\address{Division of Mathematical Sciences, School of Physical and Mathematical Sciences, Nanyang Technological University, Singapore 637371}
\email{gmpan@ntu.edu.sg}

\address{Department of Statistics and Applied Probability, National University of
 Singapore, Singapore 117546}
\email{stazw@nus.edu.sg}


\date{\today}


\maketitle

\begin{abstract}
In this paper, we study the complex Wigner matrices $M_n=\frac{1}{\sqrt{n}}W_n$ whose eigenvalues are typically in the interval $[-2,2]$. Let $\lambda_1\leq \lambda_2\cdots\leq\lambda_n$ be the ordered eigenvalues of $M_n$. Under the assumption of four matching moments with the Gaussian Unitary Ensemble(GUE), for test function $f$ 4-times continuously differentiable on an open interval including $[-2,2]$, we establish  central limit theorems for two types of partial linear statistics of the eigenvalues. The first type is defined with a threshold $u$ in the bulk of the Wigner semicircle law as $\mathcal{A}_n[f; u]=\sum_{l=1}^nf(\lambda_l)\mathbf{1}_{\{\lambda_l\leq u\}}$. And the second one is $\mathcal{B}_n[f; k]=\sum_{l=1}^{k}f(\lambda_l)$ with positive integer $k=k_n$ such that $k/n\rightarrow y\in (0,1)$ as $n$ tends to infinity. Moreover, we derive a weak convergence result for a partial sum process constructed from $\mathcal{B}_n[f; \lfloor nt\rfloor]$.
\end{abstract}
\section{Introduction}
The complex Wigner Ensemble is defined as a family of $n\times n$ random Hermitian matrices $M_n$ of the form
\begin{eqnarray*}
M_n=\frac{1}{\sqrt{n}}W_n=\frac{1}{\sqrt{n}}\{w_{jk}\}_{j,k=1}^n,
\end{eqnarray*}
in which $w_{ll}\in \mathbb{R}, 1\leq l\leq n$, $w_{jk}=\bar{w}_{kj}\in \mathbb{C},1\leq j<k\leq n$, and $\{w_{ll}, w_{jk}; 1\leq i \leq n, 1\leq j<k\leq n\}$ is a collection of independent variables such that
\begin{eqnarray*}
\mathbb{E}w_{ll}=\mathbb{E}w_{jk}=0,\quad \mathbb{E}|w_{jk}|^2=1,\quad \mathbb{E}w_{ll}^2=\sigma^2<\infty.
\end{eqnarray*}

Our basic additional assumption on the elements of $W_n$ throughout the paper is the following condition.\\\\
{\emph{Condition }}{\bf{$\mathbf{C_0}$}}:\emph{ We say that a complex Wigner matrix $M_n$ obeys Condition $\mathbf{C_0}$ if $\{w_{ll}, \mathrm{Re}w_{jk},\mathrm{Im}w_{jk}; 1\leq l \leq n, 1\leq j<k\leq n\}$ is a collection of independent variables whose distributions are all supported on at least three points, and we have the exponential decay condition on the elements in the sense that
\begin{eqnarray*}
\mathbb{P}(|w_{jk}|\geq t^C)\leq e^{-t}
\end{eqnarray*}
holds for all $t\geq C'$ with some positive constants $C,C'$ (independent of $j,k,n$ ).}\\\\

A basic example of the complex Wigner matrix satisfying {\emph{Condition}} $\mathbf{C_0}$ is drawn from the Gaussian Unitary Ensemble(GUE) whose elements are Gaussian distributed, i.e.
\begin{eqnarray*}
w_{ll}\sim N(0,1)_{\mathbb{R}},\quad 1\leq l\leq n, \quad w_{jk}\sim N(0,1)_{\mathbb{C}},\quad 1\leq j<k\leq n.
\end{eqnarray*}
Here $N(0,1)_{\mathbb{R}}$ (resp. $N(0,1)_{\mathbb{C}}$) represents the standard real (resp. complex) Gaussian distribution.

For Wigner matrix $M_n$, we denote its ordered eigenvalues as $\lambda_1(M_n)\leq \lambda_2(M_n)\leq\cdots\leq\lambda_n(M_n)$. And the empirical spectral distribution (ESD) of $M_n$ is defined by
 \begin{eqnarray*}
 F^{M_n}(x)=:\frac1n\sum_{l=1}^n\mathbf{1}_{\{\lambda_l(M_n)\leq x\}}.
 \end{eqnarray*}
 When there is no confusion, we will briefly use $\lambda_l$ and $F_n(x)$ to represent $\lambda_l(M_n)$ and $F^{M_n}(x)$ respectively.

A fundamental result is the Wigner semicircle law, which describes the global limiting behavior of eigenvalues of the Wigner ensemble: for any $\varphi\in C_b(\mathbb{R})$ (the set of bounded continuous functions in $\mathbb{R}$), one has
\begin{eqnarray}
\frac1n\sum_{l=1}^n\varphi(\lambda_l)\stackrel{\mathbb{P}}\longrightarrow \int \varphi(x)\rho_{sc}(x)dx, \label{1.1}
\end{eqnarray}
where
\begin{eqnarray*}
\rho_{sc}(x)=\frac{1}{2\pi^2}\sqrt{4-x^2}\mathbf{1}_{\{|x|\leq2\}}
\end{eqnarray*}
is the density function of the Wigner semicircle law $F_{sc}(x)$. That is to say,
the ESD $F_n(x)$ converges weakly in probability to the semicircle law $F_{sc}(x)$. We remark here (\ref{1.1}) holds under much weaker condition than $\mathbf{C_0}$ assumed in this paper, see \cite{BS} for instance. Note that (\ref{1.1}) can be viewed as a universal result corresponding to the classical law of large number (LLN) for sums of independent random variables. The quantity
\begin{eqnarray*}
\mathcal{L}_n[\varphi]=\sum_{l=1}^n\varphi(\lambda_l)
\end{eqnarray*}
is usually referred to as the global linear eigenvalue statistic ({\bf{GLES}}) of Wigner matrices with test function $\varphi$.

Once the LLN was obtained, a natural question in the probability theory is to study the fluctuation of $\mathcal{L}_n[\varphi]$ subsequently. For any smooth enough test function $\varphi$, there are a vast of results obtained on the central limit theorem for $\mathcal{L}_n[\varphi]$ under different assumptions, for instance, see \cite{BWZ}, \cite{BY}, \cite{SC}, \cite{LP}, \cite{S}, \cite{SS}. A remarkable work on this topic is due to Lytova and Pastur \cite{LP}. Particularly for GUE, Lytova and Pastur showed that for any bounded test function $\varphi$ with bounded derivative, one has
\begin{eqnarray}
\mathcal{L}_n[\varphi]-\mathbb{E}\mathcal{L}_n[\varphi]\stackrel{d}\longrightarrow N(0,V_{GUE}[\varphi]), \label{1.0}
\end{eqnarray}
where
\begin{eqnarray}
V_{GUE}[\varphi]=\frac{1}{4\pi^2}\int_{-2}^2\int_{-2}^2\left(\frac{\varphi(\lambda)-\varphi(\mu)}{\lambda-\mu}\right)^2
\frac{4-\lambda\mu}{\sqrt{4-\lambda^2}\sqrt{4-\mu^2}}d\lambda d\mu. \label{1.2}
\end{eqnarray}
(See Remark 2.1 of \cite{LP}). Moreover, an analogous result for more general Wigner matrices can be derived through the discussion in \cite{LP} for essentially $C^5$ test functions. As the reader might notice, there is no normalizing constant in the convergence in (\ref{1.0}). The reason behind this is mainly that the eigenvalues repel each other and so are more regularly distributed than independent random variables.

The main aim of this paper is to study the CLTs for two types of partial linear eigenvalue statistics ({\bf{PLES}}) in the sense that only a part of eigenvalues will be involved in the statistics.
The type 1 {\bf{PLES}} with the test function $f$ and the threshold $u\in[-2+\delta,2-\delta]$ is defined by
\begin{eqnarray*}
\mathcal{A}_n[f;u]=\sum_{l=1}^nf(\lambda_l)\mathbf{1}_{\{\lambda_l\leq u\}},
\end{eqnarray*}
which is a summation of $f(\lambda_l)$ only for $\lambda_l\leq u$. The type 2 {\bf{PLES}} with the test function $f$ and the integer $k=:k_n$ is defined by
\begin{eqnarray*}
\mathcal{B}_n[f;k]=\sum_{l=1}^kf(\lambda_l),
\end{eqnarray*}
with the constraint that $k/n\rightarrow y\in (0,1)$ as $n$ tends to infinity. Note that the type 1 {\bf{PLES}} $\mathcal{A}_n[f;u]$ is just a {\bf{GLES}} with the probably discontinuous test function $f(x)\mathbf{1}_{\{x\leq u\}}$. When $f(u)=0$, $f(x)\mathbf{1}_{\{x\leq u\}}$ is continuous but may be non differentiable at $u$. We remark here though we define the two types of {\bf{PLES}} as the sum of $f(\lambda_l)$ for the smallest eigenvalues, it will cause no intrinsic difference on all discussions throughout the paper if we define the {\bf{PLES}} with the largest eigenvalues instead of the smallest ones.

{\bf{PLES}} for a variety of matrices (deterministic or random) play relevant roles in a lot of fields. For instance, when $f(x)\equiv 1$, the type 1 {\bf{PLES}} is just the counting function of the eigenvalues up to $u$, which is a fundamental and well studied quantity in Random Matrix Theory (RMT). For the fluctuation of the counting function of Wigner matrices, we refer to \cite{CL}, \cite{G}, \cite{So}, \cite{DV} for details of this topic. For the type 2 {\bf{PLES}}, a canonical example is the sum of the $k$ largest or smallest eigenvalues, which is important in both pure and applied aspects of matrix theory. Especially, the sum of the $k$ largest eigenvalues is interesting in a lot of fields such as principal component analysis, compressed sensing and computational mathematics, see \cite{A}, \cite{BFGMR}, \cite{EAS}, \cite{OH} for instance. However, the type 2 {\bf{PLES}} are always not easy to be studied since they are concerned with the ordered eigenvalues. By the generalized Rayleigh-Ritz theorem (see Corollary 4.3.18 of \cite{HJ} for instance), one has for an $n\times n$ Hermitian matrix $A$, there exists a variational representation as
\begin{eqnarray}
\sum_{l=1}^k \lambda_l(A)=\min_{U^*U=I_k} Tr U^*AU \label{1.3}
\end{eqnarray}
with any $1\leq k\leq n$. Here $\lambda_1(A)\leq\lambda_2(A)\cdots\leq \lambda_n(A)$ are ordered eigenvalues of $A$. However, such a variational characterization is not convenient for computation and analysis. Instead, one can work with a proxy of the quantity (\ref{1.3}) by a type 1 partial sum with a threshold $u$ ``near'' $\lambda_k(A)$ as
\begin{eqnarray*}
\mathcal{A}_n[x;u]=\sum_{l=1}^n \lambda_l(A)\mathbf{1}_{\{\lambda_l(A)\leq u\}}.
\end{eqnarray*}
Such an elementary approximate technic suggests us to study the two types of {\bf{PLES}} together. As will be seen, with the aid of the so-called rigidity property of the eigenvalues proved by Erd\H{o}s, Yau and Yin \cite{EYY}, such an approximate strategy does work well for the study of the fluctuation of a type 2 {\bf{PLES}}. As mentioned above, an advantage of $\mathcal{A}_n[f;u]$ is that it is indeed a {\bf{GLES}}, though the test function is not necessarily continuous. Such a fact can help one avoid working with ordered eigenvalues.

For brevity, we set
$$f_u(x)=:(f(x)-f(u))\mathbf{1}_{\{x\leq u\}}.$$
And for $t\in [0,1]$, let $\gamma_t$ be the number that
\begin{eqnarray*}
\int_{-2}^{\gamma_t}\rho_{sc}(x)dx=t.
\end{eqnarray*}
Moreover, for some small positive number $\delta$, we set the interval
$$\mathcal{U}=(-2-\delta,2+\delta)$$
throughout the paper. We use the notation $C^{k}(\mathcal{U})$ to indicate the set of the real functions which are defined on the whole real line and $k$-times continuously differentiable on the interval $\mathcal{U}$. Our first result is the following theorem.
\begin{thm} \label{th.1} If $M_n$ is drawn from $GUE$, $u\in [-2+\delta,2-\delta]$ with some small but fixed $\delta>0$ and $k=:k_n$ such that $k/n\rightarrow y$ for some fixed constant $y\in (0,1)$ as $n$ tends to infinity, one has the following CLTs.

(i): If $f\in C^1(\mathcal{U})$ and $f(u)\neq 0$, one has
\begin{eqnarray*}
\frac{\mathcal{A}_n[f;u]-m[f;u]}{\sqrt{\frac{f^2(u)}{2\pi^2}\log n}}\stackrel d \longrightarrow N(0,1)
\end{eqnarray*}
with
\begin{eqnarray*}
m[f;u]=n\int_{-2}^uf(x)dF_{sc}(x) .
\end{eqnarray*}

(ii): If $f\in C^4(\mathcal{U})$ and $f(u)=0$, one has
\begin{eqnarray*}
\mathcal{A}_n[f;u]-m[f;u]\stackrel d \longrightarrow N(0,V_{GUE}[f_u])
\end{eqnarray*}

(iii): If $f\in C^4(\mathcal{U})$, for $\mathcal{B}_n[f;k]$ we have
\begin{eqnarray*}
\mathcal{B}_n[f;k]-m[f;\gamma_{k/n}]\stackrel d \longrightarrow N(0,V_{GUE}[f_{\gamma_y}]).
\end{eqnarray*}
\end{thm}
\begin{rem} \label{rem.1} If $f(x)$ is continuous in $\mathbb{R}$ and its magnitude grows no faster than $C_1e^{c_1x^2}$ as $|x|\rightarrow \infty$ for some constants $c_1>0$ and $C_1<\infty$, we can also replace $m[f;u]$ and $m[f;\gamma_{k/n}]$ by $\mathbb{E}\mathcal{A}_n[f;u]$ and $\mathbb{E}\mathcal{B}_n[f;k]$ respectively in the above formulas. This is a consequence of the fact that the density of $\mathbb{E}F_n(x)$ has a tail of $O(e^{-cnx^2})$. For details, see the forthcoming discussions in Section 2 and Section 3.
\end{rem}
\begin{rem} Observe that $\mathcal{A}_n[f;u]=\mathcal{L}_n[f_u]$ if $f(u)=0$. Thus (ii) implies that the differentiability condition imposed on the test function is not necessary for Lytova and Pastur's CLT, since $f_u$ may be non-differentiable at $u$. At the same time, one can also learn from (i) that a discontinuous point of the test function will indeed cause significant change on the fluctuation of {\bf{GLES}}. Note that there is still a large gap between differentiability and discontinuity. It will be interesting to investigate the relation between the smoothness of the test function and the limiting behavior of the corresponding {\bf{GLES}}.
\end{rem}

For ease of presentation, we use the notation
$$\xi^\circ=:\xi-\mathbb{E}\xi$$
for any random variable $\xi$ in the sequel. Note that (iii) of Theorem \ref{th.1} reveals the weak convergence of the random sequence $\{\mathcal{B}_n[f;k]\}$. Inspired by the classical partial sum process of i.i.d random variables, we take a step further to study the following partial sum process constructed from $\mathcal{B}_n[f;k]$ with some small but fixed $\delta>0$ as
\begin{eqnarray*}
\mathcal{S}_n[f;t]=\mathcal{B}_n[f;\lfloor nt\rfloor]+(nt-\lfloor nt\rfloor)f(\lambda_{\lfloor nt\rfloor+1}), t\in[\delta,1-\delta],
\end{eqnarray*}
which is an element of $C[\delta,1-\delta]$. Here $C[a,b]$ represents the space of continuous functions on the interval $[a,b]$ equipped with uniform topology. Then for $\mathcal{S}^\circ_n[f;t]$, we have the following weak convergence theorem in $C[\delta,1-\delta]$.
\begin{thm} \label{th.3}Suppose that $M_n$ is drawn form GUE, and $f\in C^3(\mathcal{U})$. We also assume that there exist constants $c_1>0$ and $C_1<\infty$ such that $f(x)$ is continuous in $\mathbb{R}$ and its magnitude grows no faster than $C_1e^{c_1x^2}$ as $|x|\rightarrow \infty$. Then the sequence $(\mathcal{S}^\circ_n[f;t]; t\in[\delta,1-\delta])$ is tight and converges weakly to a Gaussian process $(\mathcal{S}[f;t];t\in[\delta,1-\delta])$ with mean zero and covariance function given by
\begin{eqnarray*}
&&Cov(\mathcal{S}[f;s],\mathcal{S}[f;t])\nonumber\\
&&=\frac{1}{4\pi^2}\int_{-2}^2\int_{-2}^2\left(\frac{f_{\gamma_t}(\lambda)-f_{\gamma_t}(\mu)}{\lambda-\mu}\right)
\left(\frac{f_{\gamma_s}(\lambda)-f_{\gamma_s}(\mu)}{\lambda-\mu}\right)
\frac{4-\lambda\mu}{\sqrt{4-\lambda^2}\sqrt{4-\mu^2}}d\lambda d\mu.
\end{eqnarray*}
\end{thm}
\begin{rem} \label{rem.2} Actually, one can extend the above result to the test function $f\in C^3(U)$ without any additional condition imposed on its growth as $|x|\rightarrow \infty$ if we consider the process $\mathcal{S}[f;t]-\mathbb{E}\mathcal{S}[\bar{f}_\epsilon;t]$ instead. Here $\bar{f}_\epsilon(x)$ is a smooth truncation of $f(x)$ in the sense that $\bar{f}_\epsilon(x)=:\chi_\epsilon(x)f(x)$, where $\chi_\epsilon(x)$ is a smooth cutoff to the region $|x|\leq 2+\epsilon$ that equals $1$ for $|x|\leq 2+\epsilon/2$ with some small positive number $\epsilon<\delta$. Such an extension can be achieved easily by using the large deviation estimate of extreme eigenvalues (See Lemma \ref{lem.12} for instance). We leave the detail to the reader.
\end{rem}
\begin{rem}
It would also seem natural to consider the process constructed from the type 1 {\bf{PLES}} $\mathcal{A}_n[f;u]$ with the parameter $u\in[-2+\delta,2-\delta]$, which can be viewed as an element in $D[-2+\delta,2-\delta]$ (the metric space of functions on $[-2+\delta,2-\delta]$ with discontinuities of the first kind, equipped with Skorokhold metric.) However, we assert that for general $f$, $\mathcal{A}^\circ_n[f;u]$ cannot converge weakly to any non-trivial process for any choice of normalization. For instance, when $f(x)\equiv 1$, such a fact has been mentioned in Bai and Silverstein \cite{BS} for sample covariance matrices. The case of Wigner matrices is just analogous.
\end{rem}
The next two results are the extensions of Theorem \ref{th.1} and Theorem \ref{th.3} to more general complex Wigner matrices. In order to state our results, we use the terminology of Tao and Vu (see \cite{TV2} for example) to say that $M_n=\frac{1}{\sqrt{n}}(w_{jk})_{j,k=1}^n$ matches $M'_n=\frac{1}{\sqrt{n}}(w'_{jk})_{j,k=1}^n$ to the $\beta$-th order off the diagonal and the $\gamma$-th order on the diagonal if
$$
\mathbb{E}(w_{ll})^{\alpha}=\mathbb{E}(w'_{ll})^\alpha,\quad 0\leq\alpha\leq \gamma,\quad 1\leq l\leq n,
$$
$$
\mathbb{E}({\rm{Re}}w_{jk})^{\alpha_1}({\rm{Im}}w_{jk})^{\alpha_2}=\mathbb{E}({\rm{Re}}{w'_{jk}})^{\alpha_1}({\rm{Im}}{w'_{jk}})^{\alpha_2}, \quad 0\leq \alpha_1+\alpha_2\leq \beta, \quad 1\leq j<k\leq n,
$$
where $\alpha_1,\alpha_2$ and $\alpha$ are non-negative integers. We state our results as follows.
\begin{thm} \label{th.2} Suppose that $M_n$ is a Wigner matrix satisfying Condition $\mathbf{C_0}$ and matches to GUE to the fourth order off the diagonal and the second order on the diagonal. Then for the test function $f\in C^4(\mathcal{U})$, $(i)-(iii)$ of Theorem \ref{th.1} still hold for $M_n$.
\end{thm}
\begin{thm} \label{th.4}Suppose that $M_n$ is a Wigner matrix satisfying Condition $\mathbf{C_0}$ and matches to GUE to the fourth order off the diagonal and the second order on the diagonal. We assume $f\in C^4(\mathcal{U})$. Additionally, we assume there exist constants $K<\infty$ and $C_1<\infty$ independent of $n$, such that $f(x)$ is continuous in $\mathbb{R}$ and its magnitude grows no faster than $C_1|x|^K$ when $|x|\rightarrow \infty$. Then we also have that the process $(\mathcal{S}^\circ_n[f;t]; t\in[\delta,1-\delta])$ is tight and converges weakly to $(\mathcal{S}[f;t];t\in[\delta,1-\delta])$.
\end{thm}
\begin{rem} Similar to Remark \ref{rem.2}, one can also extend the above result to $f\in C^4(\mathcal{U})$ if one considers the process $\mathcal{S}[f;t]-\mathbb{E}\mathcal{S}[\bar{f}_\epsilon;t]$ instead. Moreover, it is likely that one can extend the above result to $t\in [\delta,1]$ with further discussion on the edge of the spectrum. However, the current issue relies on some crucial estimates proved only for the bulk case, such as Lemma \ref{lem.2}. So we do not pursue this direction here.
\end{rem}
From now on, we will use the notation $C,C_1,C'$ and $L$ to denote some $n$-independent positive constants whose values may differ from line to line. And throughout the paper, we say an event $E$ holds with high probability if
\begin{eqnarray*}
P(E)\geq 1-n^{-c}
\end{eqnarray*}
with some constant $c$ and with overwhelming probability if
\begin{eqnarray*}
P(E)\geq 1-n^{-C}
\end{eqnarray*}
for any constant $C>0$.

Our paper is organized as follows. In Section 2, we provide some basic tools and preliminaries of the whole paper. And Section 3 is devoted to the proof of the CLTs for two types of {\bf{PLES}} for GUE, i.e. Theorem \ref{th.1}. Then in Section 4, we prove Theorem \ref{th.3}, whose proof is heavily based on the discussion in Section 3. In Section 5, we prove a comparison theorem for linear eigenvalue statistics, see Theorem \ref{th.3.1}. And as an application, we use our comparison theorem to extend Theorem \ref{th.1} to general complex Wigner matrices, i.e. Theorem \ref{th.2}. Also with the aid of the comparison theorem, we prove Theorem \ref{th.4} in Section 6. Some necessary known results are stated in the Appendix.
\section{Preliminaries}
In this section, we will provide some basic notions and tools necessary for our proof in the sequel. Totally speaking, our strategy is to prove the results for GUE first and then extend them to general Wigner matrices by some comparison procedure. Such a strategy is quite fundamental in RMT. Thus the basic tools presented in this section consist of two separated parts. The first part is particularly for GUE, and the second part will be mainly contributed to our comparison procedure.

Using GUE as our starting point is mainly because its explicit formula of the joint probability density (j.p.d.) for the eigenvalues has a determinantal structure, which is friendly with analysis. By making use of the j.p.d., a vast of central issues in RMT can be solved explicitly for GUE. We refer to the books of Deift \cite{De} and Mehta \cite{Me} for comprehensive surveys in this aspect.

If $M_n=\frac{1}{\sqrt{n}}W_n$ is drawn from GUE, then the joint distribution of non ordered eigenvalues of $W_n$ has the following j.p.d.
\begin{eqnarray*}
\rho_{n,n}(x_1,\cdots,x_n)=\frac{1}{n!}\prod_{1\leq j<k\leq n}|x_j-x_k|^2e^{-\sum_{j=1}^nx_j^2/2}.
\end{eqnarray*}
For the point process $x_1,\cdots,x_n$, the $k$-point correlation function $\rho_{k,n}$ has the well known determinantal structure
\begin{eqnarray*}
\rho_{k,n}(x_1,\cdots,x_k)
&=:&\frac{n!}{(n-k)!}\int_{\mathbb{R}^{n-k}}\rho_{n,n}(x_1,\cdots,x_n)dx_{k+1}\cdots dx_n\nonumber\\
&=&\det(K_n(x_l,x_j))_{l,j=1}^k.
\end{eqnarray*}
Here $K_n(x,y)$ is the kernel function given by
\begin{eqnarray*}
K_n(x,y)=\sum_{l=0}^{n-1}H_l(x)H_l(y)e^{-\frac14(x^2+y^2)}=\sum_{i=0}^{n-1}\psi_l(x)\psi_l(y),
\end{eqnarray*}
where $H_l(x)$ is the $l$-th orthonormalized Hermite polynomial w.r.t. the weight function $e^{-x^2/2}$ and $\psi_l(x)$ is the corresponding oscillator wave function. By the famous Christoffel-Darboux formula, one has for $x\neq y$
\begin{eqnarray*}
K_n(x,y)=\sqrt{n}\frac{\psi_n(x)\psi_{n-1}(y)-\psi_{n-1}(x)\psi_n(y)}{x-y}
\end{eqnarray*}
and for $x=y$ by l'H\^{o}pital's rule,
\begin{eqnarray}
K_n(x,x)=n\psi_{n-1}^2(x)-\sqrt{n(n-1)}\psi_{n-2}(x)\psi_n(x).\label{2.27}
\end{eqnarray}
Using the notation
\begin{eqnarray*}
\mathcal{K}_n(x,y)=\sqrt{n}K_n(\sqrt{n}x,\sqrt{n}y),
\end{eqnarray*}
one has the following explicit formulas of expectation and variance of $\mathcal{L}_n[\varphi]$ for GUE.
\begin{eqnarray}
\mathbb{E}\mathcal{L}_n[\varphi]=\int_{\mathbb{R}}\varphi(x)\mathcal{K}_n(x,x)dx \label{2.25}
\end{eqnarray}
and
\begin{eqnarray}
Var\mathcal{L}_n[\varphi]=\frac12\int_{\mathbb{R}}\int_{\mathbb{R}} (\varphi(x)-\varphi(y))^2\mathcal{K}^2_n(x,y)dxdy. \label{2.26}
\end{eqnarray}
See (4) of \cite{So} and (2.27) of \cite{Pa} for instance. Note that from (\ref{2.25}), $\mathcal{K}_n(x,x)/n$ is just the density function of $\mathbb{E}F_n(x)$. For our purpose, we state below some properties for the kernel function $\mathcal{K}_n(x,y)$.
Firstly, by definition and the Cauchy-Schwarz inequality,
\begin{eqnarray}
|\mathcal{K}_n(x,y)|^2\leq \mathcal{K}_n(x,x)\mathcal{K}_n(y,y) \label{2.28}
\end{eqnarray}
Moreover, by (\ref{2.27}), one has
\begin{eqnarray*}
\mathcal{K}_n(x,x)=\sqrt{n}(n\psi_{n-1}^2(\sqrt{n}x)-\sqrt{n(n-1)}\psi_{n-2}(\sqrt{n}x)\psi_n(\sqrt{n}x)).
\end{eqnarray*}
By adjusting the scale in the setting of \cite{G}, one can see that when $|x|\geq 2+\varepsilon$ for some small $\varepsilon>0$,
\begin{eqnarray}
\psi_n(\sqrt{n}x)=\mathcal{O}(n^{-1/4}e^{-\frac{n}{4}F(x)}),\label{2.29}
\end{eqnarray}
where
\begin{eqnarray*}
F(x)=|\int_x^2\sqrt{|4-y^2|}dy|.
\end{eqnarray*}
See Section 4 of \cite{G} for reference. Consequently,
\begin{eqnarray}
\frac{\mathcal{K}_n(x,x)}{n}=\mathcal{O}(e^{-cnx^2}) \label{2.33}
\end{eqnarray}
for $|x|\geq 2+\varepsilon$. Such a fact has been mentioned in Remark \ref{rem.1}.
By (\ref{2.25})-(\ref{2.33}), one can see that when $\varphi(x)$ is continuous and $|\varphi(x)|$ grows more slowly than $C_1e^{c_1x^2}$ as $|x|\rightarrow \infty$ for some $c_1>0$ and $C_1<\infty$, we have for sufficiently large $n$
\begin{eqnarray}
\mathbb{E}\mathcal{L}_n[\varphi]=\int_{\mathcal{U}}\varphi(x)\mathcal{K}_n(x,x)dx+\mathcal{O}(e^{-cn}) \label{2.30}
\end{eqnarray}
and
\begin{eqnarray}
Var\mathcal{L}_n[\varphi]=\frac12\int_{\mathcal{U}}\int_{\mathcal{U}} (\varphi(x)-\varphi(y))^2\mathcal{K}^2_n(x,y)dxdy+\mathcal{O}(e^{-cn}) \label{2.31}
\end{eqnarray}
with some positive constant $c$ depending only on $\varphi$ and $\delta$. The above formulas will be frequently used in our proof for the GUE case in Sections 3 and 4.

However, for general Wigner matrices, the explicit formula for the joint distribution of the eigenvalues is obviously not available. A classical strategy in probability theory is the so-called Lindeberg method to replace a non-Gaussian variable by a Gaussian one at each step, and to study the stability of the concerned quantity under such a swapping procedure. A successful use of Lindeberg method to RMT in the recent work of Tao and Vu \cite{TV1} helped to extend a lot of results on local eigenvalue statistics from GUE to general Wigner matrices. However, Tao and Vu's strategy in \cite{TV1} requires a detailed analysis on the spectral dynamics of the matrices in the sense that the accurate estimates of the derivatives of the eigenvalues w.r.t the matrix elements are needed.

Later on, Erd\H{o}s, Yau and Yin proposed another swapping strategy to derive the bulk universality of local statistics in \cite{EYY1}. They studied the stability of the Green function instead of eigenvalues under every swapping step. Such a strategy is based on the elementary resolvent expansion formula (see (\ref{2.34})) and turns out to be relatively simpler for certain problems. Very recently, Tao and Vu used a similar swapping strategy on the Green function to derive the CLT for the log-determinant and a sharp concentration of counting functions for Wigner matrices, see \cite{TV2} and \cite{TV3}. Note that the objects in \cite{TV2} and \cite{TV3} are just two examples of {\bf{GLES}} with discontinuous test functions (logarithmic and indicator function respectively). It will be clear that one major technical difficulty in our problem is to derive a CLT for {\bf{GLES}} with the test function continuous but maybe non-differentiable at a few points. Such an ill behaviour in smoothness leads us to pursue the idea in \cite{TV2} and \cite{TV3} to study the {\bf{GLES}} with a class of non-smooth test functions. In Section 5, we will establish a comparison theorem for linear eigenvalue statistics, based on the Helffer-Sj\"{o}strand formula and resolvent expansion. For this purpose, we state some related notions and tools in the remaining part of this section.

The Stieltjes transform of a probability measure $\mu$ can be defined for all complex number $z\in \mathbb{C}\setminus\mathbb{R}$ as
\begin{eqnarray*}
s^{\mu}(z)= \int \frac{1}{x-z}\mu(dx).
\end{eqnarray*}
Thus for the ESD $F_n(x)$ we have
\begin{eqnarray*}
s^{M_n}(z)=:s^{F_n}(z)=\frac1n\sum_{l=1}^n\frac{1}{\lambda_l(M_n)-z}=\frac1nTr(M_n-zI_n)^{-1}.
\end{eqnarray*}
And we also denote the resolvent of $M_n$ by
\begin{eqnarray*}
R^{M_n}(z)=(M_n-zI_n)^{-1},
\end{eqnarray*}
Thus $s^{M_n}(z)=\frac1nTrR^{M_n}(z)$. When there is no confusion, we will simplify the symbols $s^{M_n}(z), R^{M_n}(z)$ by $s_n(z), R_n(z)$.

Using the terminology of \cite{TV2}, we say a matrix $V$ is an elementary matrix if it has one of the following forms
 \begin{eqnarray*}
 V=e_je_j^*, e_je_k^*+e_ke_j^*, ie_je_k^*-ie_ke_j^*
 \end{eqnarray*}
 with $1\leq j\neq k\leq n$. Here $e_1,\cdots, e_n$ is the standard basis of $\mathbb{C}^n$. Let $M_0$ be an $n\times n$ Hermitian matrix, and set $M_t=M_0+\frac{1}{\sqrt{n}}tV$. Correspondingly, we denote the resolvent and Stieltjes transform of $M_t$ by $R_t(z)$ and $s_t(z)$ respectively for some complex number $z=x+iy$ with $y\neq 0$. When there is no confusion, we will simplify the notation $R_t(z),s_t(z)$ by $R_t,s_t$. The notation $||A||_{(\infty,1)}$ for a matrix $A=(a_{jk})_{j,k=1}^n$ means its $l^{1}\rightarrow l^\infty$ operator norm in the sense that
 \begin{eqnarray*}
 ||A||_{(\infty,1)}=\sup_{1\leq j,k\leq n}|a_{jk}|.
 \end{eqnarray*}
 We conclude this section by the following crucial Taylor expansion for $s_t$ provided by Tao and Vu.

 \begin{lem} \label{lem.3}\emph{(Proposition 13, \cite{TV2})}Suppose that $x\in\mathbb{R}$, $y>0$ and $t\in\mathbb{R}$. If
 \begin{eqnarray}
 |t|||R_0||_{(\infty,1)}=o(\sqrt{n}),\label{3.2}
 \end{eqnarray}
 one has for fixed integer $k\geq 0$,
 \begin{eqnarray*}
 s_t=s_0+\sum_{j=1}^kn^{-j/2}c_jt^j+\mathcal{O}\left(n^{-(k+1)/2}|t|^{k+1}||R_0||^{k+1}_{(\infty,1)}\min (||R_0||_{(\infty,1)},\frac{1}{ny})\right)
 \end{eqnarray*}
 where the coefficients $c_j$ are independent of $t$ and obey the bounds
 \begin{eqnarray*}
 c_j=\mathcal{O}\left(||R_0||^j_{(\infty,1)}\min \{||R_0||_{(\infty,1)},\frac{1}{ny}\}\right)
 \end{eqnarray*}
 for all $1\leq j\leq k$.
 \end{lem}

 Lemma \ref{lem.3} is a consequence of the elementary resolvent expansion formula
 \begin{eqnarray}
 R_t=R_0+\sum_{j=1}^k\left(-\frac{t}{\sqrt{n}}\right)^j(R_0V)^jR_0+\left(-\frac{t}{\sqrt{n}}\right)^{k+1}(R_0V)^{k+1}R_t. \label{2.34}
 \end{eqnarray}
 We refer to \cite{TV2} for the details of the proof.
\section{CLTs For Gaussian Case}
First, we truncate the test function so that it is compactly supported and show that such a modification does not alter our results. Set the interval $\mathcal{U}_\epsilon=[-2-\epsilon,2+\epsilon]$ with a small constant $\epsilon<\delta$. For test function $f(x)$, we define the truncated function $\bar{f}_\epsilon(x)=:\chi_\epsilon(x)f(x)$. Here $\chi_\epsilon(x)$ is a smooth cutoff to the region $|x|\leq 2+\epsilon$ that equals $1$ for $|x|\leq 2+\epsilon/2$. It follows from Lemma \ref{lem.12} in Appendix that
\begin{eqnarray*}
\mathbb{P}(\mathcal{A}_n[f;u]\neq \mathcal{A}_n[\bar{f}_\epsilon;u])\leq \mathbb{P}(\max_{1\leq l\leq n}|\lambda_l|\geq 2+\epsilon/2)\rightarrow 0
\end{eqnarray*}
as $n$ goes to infinity. Consequently, without loss of generality, we can thus always assume that $f$ is compactly supported on the interval $\mathcal{U}_{\epsilon}\subset \mathcal{U}$ in this Section.

To prove Theorem \ref{th.1}, we start with the type 1 {\bf{PLES}}. We do the decomposition as follows
\begin{eqnarray*}
\mathcal{A}_n[f;u]=\sum_{l=1}^n(f(\lambda_l)-f(u))\mathbf{1}_{\{\lambda_l\leq u\}}+f(u)N_n(-\infty,u].
\end{eqnarray*}
With the notation defined above, we have
\begin{eqnarray*}
\mathcal{A}_n[f;u]=\sum_{l=1}^n f_u(\lambda_l)+f(u)N_n(-\infty,u]=\mathcal{L}_n[f_u]+f(u)N_n(-\infty,u].
\end{eqnarray*}
Observe that $f_u(x)$ is a continuous function with only one possibly non-differentiable point $u$. In order to apply the approach in \cite{LP} to treat such a test function, we smooth $f_u(x)$ in a tiny interval including $u$. Set the interval
$$I_n(u)=[u-n^{-1/2-c},u+n^{-1/2-c}]:=[a_n(u),b_n(u)],$$
where $c$ is a small positive constant. Define the smooth modification of $f_u(x)$ by
\begin{eqnarray*}
g_u(x)=:g_u(n,x)=(f(x)-f(u))\chi_u(n,x)
\end{eqnarray*}
where $\chi_u(n,x)$ is an $n$-dependent smooth cutoff to the region $x\in(-\infty,b_n(u)]$ that equals $1$ for $x\in (-\infty, a_n(u))$, and has the property
 \begin{eqnarray}
 |\frac{d^{k}}{dx^{k}}\chi_u(n,x)|=\mathcal{O}(n^{k(1/2+c)}). \label{2.17}
 \end{eqnarray}
Consequently, one has for $f\in C^m(\mathcal{U})$,
\begin{eqnarray}
|\frac{d^{k}}{dx^{k}}g_u(x)|=\mathcal{O}(n^{(k-1)(1/2+c)}), \quad x\in I_n(u),\quad k=1,\cdots,m. \label{2.35}
\end{eqnarray}
Using Lemma \ref{lem.5.1} in Appendix one has for some positive constants $C$
\begin{eqnarray*}
N_n(I)\leq Cn|I|
\end{eqnarray*}
with overwhelming probability for any interval $I$ with length $|I|\geq n^{-1+c}$. Together with the trivial fact that
\begin{eqnarray*}
\sup_{x}|f_u(x)-g_u(x)|\leq n^{-1/2-c},
\end{eqnarray*}
we obtain
\begin{eqnarray}
|\mathcal{L}_n[f_u]-\mathcal{L}_n[g_u]|\leq \sup_{x}|f_u(x)-g_u(x)| N_n(I_n)=o(1) \label{2.8}
\end{eqnarray}
with overwhelming probability.
Consequently, we have
\begin{eqnarray}
\mathcal{A}_n[f;u]=\mathcal{L}_n[g_u]+f(u)N_n(-\infty,u]+o(1) \label{2.1}
\end{eqnarray}
holding with overwhelming probability. Furthermore, we also have
\begin{eqnarray}
\mathbb{E}\mathcal{L}_n[f_u]=\mathbb{E}\mathcal{L}_n[g_u]+o(1). \label{1.4}
\end{eqnarray}

As we have mentioned in (\ref{1.0}), for $n$-independent test function $\varphi\in C^1_b(\mathbb{R})$ with bounded derivative, Lytova and Pastur have proved the CLT. Unfortunately, here our modified test function $g_u(x)$ is $n$-dependent. Thus we can not use Lytova and Pastur's result directly. However, we will show that a slight adjustment of Lytova and Pastur's issue can still lead to the limiting behavior of $\mathcal{L}_n[g_u]$. We formulate our conclusion as the following lemma.
\begin{lem} \label{lem.6} If $M_n$ is drawn from GUE, then for $f\in C^3(\mathcal{U})$ and compactly supported on $\mathcal{U}_\epsilon$, one has
\begin{eqnarray*}
\mathcal{L}_n[g_u]-\mathbb{E}\mathcal{L}_n[g_u]\stackrel d\longrightarrow N(0,V_{GUE}[f_u]).
\end{eqnarray*}
\end{lem}
\begin{proof}
Since we follow the argument of Lytova and Pastur in \cite{LP} with only some minor changes, we sketch the proof below. Firstly, we present here some notation and known results laid out in \cite{LP}. Let
\begin{eqnarray*}
U(t)=e^{itM_n},&&u_n(t)=TrU(t),\nonumber\\
e_n(x)=e^{ix \mathcal{L}_n^\circ[g_u]}, &&Y_n(x,t)=\mathbb{E}\{u_n^\circ(t)e_n(x)\}.
\end{eqnarray*}
The basic idea of \cite{LP} is to use the characteristic function to derive a CLT. Set
\begin{eqnarray*}
Z_n(x)=:\mathbb{E}e^{ix\mathcal{L}_n^\circ[g_u]}=\mathbb{E}e_n(x).
\end{eqnarray*}
Thus it suffices to show that for any $x\in\mathbb{R}$
\begin{eqnarray*}
\lim_{n\rightarrow \infty}Z_n(x)=Z(x),
\end{eqnarray*}
where
\begin{eqnarray*}
Z(x)=\exp\{-x^2V_{GUE}[f_u]/2\}.
\end{eqnarray*}
Note the relations
\begin{eqnarray*}
Z(x)=1-V_{GUE}[f_u]\int_0^xyZ(y)dy
\end{eqnarray*}
and
\begin{eqnarray*}
Z_n(x)=1+\int_0^xZ'_n(y)dy.
\end{eqnarray*}
Using the Cauchy-Schwarz inequality and Proposition 2.4 of \cite{LP}, we see that
\begin{eqnarray*}
|Z'_n(x)|=|i\mathbb{E}\{\mathcal{L}_n^\circ[g_u]e^{ix\mathcal{L}_n^\circ[g_u]}\}|\leq \sqrt{2}\sup_{\lambda}|g'_u(\lambda)|=\mathcal{O}(1)
\end{eqnarray*}
Thus by the dominated convergence theorem, it suffices to verify that any convergent subsequences $\{Z_{n_j}\}$ and $\{Z'_{n_j}\}$ satisfy
\begin{eqnarray}
\lim_{n_j\rightarrow\infty}Z_{n_j}(x)=Z(x),\quad \lim_{n_j\rightarrow\infty}Z'_{n_j}(x)=-xV_{GUE}[f_u]Z(x). \label{2.15}
\end{eqnarray}
If we denote the Fourier transform of a function $\varphi$ by
\begin{eqnarray*}
\widehat{\varphi}(t)=\frac{1}{2\pi}\int e^{-it\lambda}\varphi(\lambda)d\lambda,
\end{eqnarray*}
 we have
 \begin{eqnarray}
 Z'_n(x)=i\mathbb{E}\{\mathcal{L}_n^\circ[g_u]e^{ix\mathcal{L}_n^\circ[g_u]}\}=i\int\widehat{g_u}(t)Y_n(x,t)dt. \label{2.32}
 \end{eqnarray}
 As shown in \cite{LP}, to prove (\ref{2.15}) one needs to prove that the sequence $\{Y_n\}$ is bounded and equicontinuous on any compact subset of $\{t\geq 0,x\in\mathbb{R}\}$ (the case of $t\leq 0$ is analogous), and every uniformly convergent on the set subsequence has the same limit $Y$. The proofs for boundness and equicontinuity are really the same as those in \cite{LP}. In fact, by the estimates in \cite{LP}, one has
 \begin{eqnarray}
 &&Var \{u_n(t)\}\leq 2t^2,\quad Var\{u'_n(t)\}\leq2(1+2t^2),\quad |Y_n(x,t)|\leq \sqrt{2}|t|,\label{2.11}\\
 &&\left|\frac{\partial}{\partial t}Y_n(x,t)\right|\leq \sqrt{2}(1+2t^2)^{1/2},\quad \left|\frac{\partial}{\partial x}Y_n(x,t)\right|\leq 2|t|\sup_{\lambda\in \mathbb{R}}|g'_u(\lambda)|\leq C t\label{2.12}.
 \end{eqnarray}
Thus the main task is to show that any uniformly convergent subsequence of $\{Y_n\}$ has the same limit $Y$, and determine the limit. A detailed estimation is presented for Gaussian Orthogonal Ensemble (GOE) in \cite{LP}. It is easy to adjust the discussion to GUE case. Applying the calculation procedure of \cite{LP} to GUE one can get
\begin{eqnarray*}
Y_n(x,t)&=&-n^{-1}\int_0^tdt_1\int_0^{t_1}\mathbb{E}\{u_n(t_2-t_1)u_n(t_2)e_n^\circ(x)\}dt_2\nonumber\\
&&-x\int_0^t\mathbb{E}\{e_n(x)n^{-1}Tr U(t_1)g'_u(M)\}dt_1.
\end{eqnarray*}
The above equation is just analogous to the corresponding one of the GOE case stated in \cite{LP}. Such a representation is a consequence of the integration by parts formula of the Gaussian variables. We refer to \cite{LP} for detail. The above equation can be rewritten as
\begin{eqnarray}
&&Y_n(x,t)+\int_0^tdt_1\int_0^{t_1}\bar{v}_n(t_1-t_2)Y_n(x,t_2)dt_2\nonumber\\
&&=xZ_n(x)A_n(t)+r_n(x,t), \label{2.14}
\end{eqnarray}
where
\begin{eqnarray*}
&&\bar{v}_n(t)=n^{-1}\mathbb{E}u_n(t), \nonumber\\
&& A_n(t)=-\int_0^t\mathbb{E}\{n^{-1}Tr U(t_1)g'_u(M)\}dt_1
\end{eqnarray*}
and
\begin{eqnarray}
r_n(x,t)
&=& -n^{-1}\int_0^{t}dt_1\int_0^{t_1}\mathbb{E}\{u_n^\circ(t_1-t_2)u_n^\circ(t_2)e_n^\circ(x)\}dt_2\nonumber\\
&&-ixn^{-1}\int_0^tdt_1\int t_2\widehat{g_u}(t_2)\mathbb{E}\{u_n(t_1+t_2)e^\circ_n(x)\}dt_2.\label{2.12}
\end{eqnarray}
By the boundness of $Y_n(x,t)$, the first inequality of (\ref{2.11}) and the Cauchy-Schwarz inequality, one immediately gets that the first term in the expression of $r_n(x,t)$ is negligible. Now we show that the second term is also $o(1)$ uniformly in any compact subset of $\{t\geq 0, x\in\mathbb{R}\}$. It suffices to prove
\begin{eqnarray}
\int_0^tdt_1\int t_2\widehat{g_u}(t_2)\mathbb{E}\{u_n(t_1+t_2)e^\circ_n(x)\}dt_2=o(n).\label{2.13}
\end{eqnarray}
In view of (\ref{2.11}),
\begin{eqnarray*}
&&|\mathbb{E}\{u_n(t_1+t_2)e^\circ_n(x)\}|=|\mathbb{E}\{u^\circ_n(t_1+t_2)e_n(x)\}|\nonumber\\
&&\leq Var^{1/2}\{u_n(t_1+t_2)\}\leq\sqrt{2}|t_1+t_2|.
\end{eqnarray*}
Hence it suffices to show that
\begin{eqnarray*}
\int(1+|t|^2)|\widehat{g_u}(t)|dt=o(n).
\end{eqnarray*}
Apparently, we can show for some positive constant $C$
\begin{eqnarray*}
\int_{|t|\geq C}(1+|t|^2)|\widehat{g_u}(t)|dt=o(n)
\end{eqnarray*}
instead.
Note that
\begin{eqnarray*}
&&\int_{|t|\geq C}(1+|t|^2)|\widehat{g_u}(t)|dt=\int_{|t|\geq C}\frac{1+t^2}{|t|^3}|\widehat{g^{(3)}_u}(t)|dt\nonumber\\
&&\leq \left(\int_{|t|\geq C}\left(\frac{1+t^2}{|t|^3}\right)^2dt\right)^{1/2}\left(\int|\widehat{g^{(3)}_u}(t)|^2dt\right)^{1/2}\nonumber\\
&&\leq C\left(\int_{\mathcal{U}_{\epsilon}}|g^{(3)}_u(x)|^2dx\right)^{1/2}\nonumber\\
&& \leq Cn^{\frac34+\frac32 c},
\end{eqnarray*}
where $g_u^{(3)}(x)=\frac{d^3}{dx^3}g_u(x)$ and $\widehat{g^{(3)}_u}(t)$ is its Fourier transform. In the above second inequality, we have used the Plancherel's Theorem, and in the last step, we used the bound (\ref{2.35}). Thus if we choose $c$ sufficiently small, we can get (\ref{2.13}). Consequently, we have
\begin{eqnarray*}
\lim_{n\rightarrow\infty}r_n(x,t)=0
\end{eqnarray*}
uniformly on any compact subset of $\{t\geq 0,x\in \mathbb{R}\}$.

Moreover, it is not difficult to derive that on any finite interval of $\mathbb{R}$, $\{\bar{v}_n\}$ and $\{A_n\}$ converge uniformly to
\begin{eqnarray}
v(t)=\frac{1}{2\pi}\int_{-2}^2e^{it\lambda}\sqrt{4-\lambda^2}d\lambda
\end{eqnarray}
and
\begin{eqnarray}
A(t)=-\frac{1}{2\pi}\int_0^t\int_{-2}^ue^{it_1\lambda}f'(\lambda)\sqrt{4-\lambda^2}d\lambda.
\end{eqnarray}
In fact, the convergence of $\bar{v}_n$ is a direct consequence of (\ref{1.1}). For the convergence of $A_n(t)$, one can use the convergence rate for ESD of GUE as
 \begin{eqnarray}
 \sup_{x}|\mathbb{E}F_n(x)-F_{sc}(x)|\leq Cn^{-1}, \label{2.16}
 \end{eqnarray}
which was proved by G\"{o}tze and Tikhomirov in \cite{GT}. Note that
\begin{eqnarray*}
A_n(t)=-\int_0^tdt_1\int e^{it_1x}g'_u(x)d\mathbb{E}F_n(x).
\end{eqnarray*}
Using (\ref{2.16}), by integration by parts, one can easily get
\begin{eqnarray*}
A_n(t)=-\int_0^tdt_1\int e^{it_1x}g'_u(x)dF_{sc}(x)+\mathcal{O}(n^{-1}).
\end{eqnarray*}
Then it is easy to see the right hand side of the above equation tends to $A(t)$ as $n$ goes to infinity.

Then by a routine analysis on the limiting equation of (\ref{2.14}) as that in \cite{LP}, one can get that $Y_n(x,t)$ converges to
 \begin{eqnarray*}
Y(x,t)=\frac{ixZ(x)}{2\pi^2}\int_{-2}^{u}\int_{-2}^2\frac{\sqrt{4-\lambda^2}}{\sqrt{4-\mu^2}} \frac{e^{it\lambda}-e^{it\mu}}{\lambda-\mu}f'(\lambda)d\lambda d\mu
\end{eqnarray*}
uniformly on any compact subset of $\{t\geq 0, x\in \mathbb{R}\}$

 Note that $\hat{g}_u(t)$ converges to $\hat{f}_u(t)$ uniformly in $t$. Thus by (\ref{2.32}) one can get for every convergence subsequence $\{Z_{n_l}\}_{l\geq 1}$ there exists
\begin{eqnarray*}
\lim_{n_l\rightarrow \infty}Z'_{n_l}(x)=-\frac{xZ(x)}{2\pi^2}\int_{-2}^{u}\int_{-2}^2\frac{\sqrt{4-\lambda^2}}{\sqrt{4-\mu^2}} \frac{f_u(\lambda)-f_u(\mu)}{\lambda-\mu}f'(\lambda)d\lambda d\mu.
\end{eqnarray*}
By the fact that for $\lambda\in (-\infty, u)$
\begin{eqnarray*}
f'(\lambda)(f_u(\lambda)-f_u(\mu))=\frac12\frac{\partial}{\partial \lambda}(f_u(\lambda)-f_u(\mu))^2,
\end{eqnarray*}
we can use integration by parts to get
\begin{eqnarray*}
\lim_{n_l\rightarrow \infty}Z'_{n_l}(x)&=&-\frac{xZ(x)}{4\pi^2}\int_{-2}^{u}\int_{-2}^2\left(\frac{f_u(\lambda)-f_u(\mu)}{\lambda-\mu}\right)^2\frac{(4-\lambda\mu)}
{\sqrt{4-\lambda^2}\sqrt{4-\mu^2}}d\lambda d\mu\nonumber\\
&&-\frac{xZ(x)}{4\pi^2}\int_{-2}^2f_u^2(\mu)\frac{\sqrt{4-u^2}}{u-\mu}\frac{1}{\sqrt{4-\mu^2}}d\mu\nonumber\\
&&=-\frac{xZ(x)}{4\pi^2}\int_{-2}^2\int_{-2}^{2}\left(\frac{f_u(\lambda)-f_u(\mu)}{\lambda-\mu}\right)^2\frac{(4-\lambda\mu)}
{\sqrt{4-\lambda^2}\sqrt{4-\mu^2}}d\lambda d\mu.
\end{eqnarray*}
Thus we conclude the proof of Lemma \ref{lem.6}.
\end{proof}
With the aid of Lemma \ref{lem.6}, we can now prove Theorem \ref{th.1}.

\begin{proof}[Proof of Theorem \ref{th.1}]
We begin with the CLT for the counting function of eigenvalues of complex Wigner matrices, whose proof can be found in the recent work of Dallaporta and Vu \cite{DV}.
\begin{lem} [\cite{DV}]  \label{lem.11}If $M_n$ is a complex Wigner matrix satisfying Condition $\mathbf{C_0}$ and matches to GUE to the fourth order off the diagonal and the second order on the diagonal, one has
\begin{eqnarray}
 \frac{N_n(-\infty,u]-\mathbb{E}N_n(-\infty,u]}{\sqrt{\frac{1}{2\pi^2}\log n}}\stackrel{d}\longrightarrow N(0,1),\label{2.4}
 \end{eqnarray}
 where
 \begin{eqnarray}
\mathbb{E}N_n(-\infty,u]=n\frac{1}{2\pi}\int_{-2}^u\sqrt{4-x^2}dx+o(1).\label{2.20}
\end{eqnarray}
\end{lem}

For $f\in C^1(\mathbb{R})$, we also have $g_u\in C^1(\mathbb{R})$. By Proposition 2.4 of \cite{LP} , we have for GUE
\begin{eqnarray*}
Var\mathcal{L}_n[g_u]\leq 2\left(\sup_{x\in \mathbb{R}}|g'_u(x)|\right)^2\leq C.
\end{eqnarray*}
Clearly,
\begin{eqnarray}
\frac{\mathcal{L}_n^\circ[g_u]}{\sqrt{\log n}}\stackrel {P}\longrightarrow 0\label{2.5}
\end{eqnarray}
as $n$ goes to infinity.
Combining (\ref{2.1}), (\ref{1.4}), (\ref{2.4}) and (\ref{2.5}), we immediately get that when $f\in C^1(\mathbb{R})$ with bounded derivative and $f(u)\neq 0$,
\begin{eqnarray*}
\frac{\mathcal{A}_n[f;u]-\mathbb{E}\mathcal{L}_n[f_u]-f(u)\mathbb{E}N_n(-\infty,u]}{\sqrt{\frac{f^2(u)}{2\pi^2}\log n}}\stackrel d\longrightarrow N(0,1).
\end{eqnarray*}
To prove (i) of Theorem \ref{th.1}, it remains to show for $f\in C^1(\mathcal{U})$ compactly supported on $\mathcal{U}_\epsilon$,
\begin{eqnarray}
\mathbb{E}\mathcal{L}_n[f_u]=n\frac{1}{2\pi}\int_{-2}^uf_u(x)\sqrt{4-x^2}dx+\mathcal{O}(1). \label{2.6}
\end{eqnarray}
Observe that by (\ref{2.16}),
\begin{eqnarray*}
|\mathbb{E}\mathcal{L}_n[f_u]-n\int_{-2}^u f_u(x) dF_{sc}(x)|
&=&n|\int_{-2-\epsilon}^u f_u(x)d\mathbb{E}F_n(x)-\int_{-2}^u f_u(x) dF_{sc}(x)|\nonumber\\
&=&n|\int_{-2-\epsilon}^u f'_u(x)(\mathbb{E}F_n(x)-F_{sc}(x))dx|\nonumber\\
&\leq&n\int_{-2-\epsilon}^u |f'_u(x)||\mathbb{E}F_n(x)-F_{sc}(x)|dx\nonumber\\
&=&\mathcal{O}(1),
\end{eqnarray*}
which implies (\ref{2.6}). Thus we complete the proof of (i).

Now we turn to the case where $f(u)=0$. Since $f\in C^4(\mathcal{U})$, by Lemma \ref{lem.6}, together with (\ref{2.1}) and (\ref{1.4}) we can easily obtain
\begin{eqnarray*}
\frac{\mathcal{A}_n[f;u]-\mathbb{E}\mathcal{L}_n[f_u]}{\sqrt{V_{GUE}[f_u]}}\stackrel d\longrightarrow N(0,1).
\end{eqnarray*}
Thus to prove (ii) of Theorem \ref{th.1}, it suffices to show for $f\in C^4(\mathcal{U})$ and compactly supported on $\mathcal{U}_\epsilon$ the more accurate estimate
\begin{eqnarray}
\mathbb{E}\mathcal{L}_n[f_u]=n\frac{1}{2\pi}\int_{-2}^uf_u(x)\sqrt{4-x^2}dx+o(1). \label{2.7}
\end{eqnarray}
To show (\ref{2.7}), we define two smooth cutoff functions $\chi_1(x)$ and $\chi_2(x)$. Let $\chi_1(x)$ be a smooth cutoff function which is equal to $1$ for $x\geq -1+u/2$ and $0$ for $x\leq -3/2+u/4$, such that $\chi_1^{(k)}(x)\leq C$ holds for some positive constant $C$ and $k=0,\cdots,4$. Let $\chi_2(x)=1-\chi_1(x)$. Now we decompose $f_u(x)$ as
\begin{eqnarray*}
f_u(x)=\chi_1(x)f_u(x)+\chi_2(x)f_u(x).
\end{eqnarray*}
Observe that $\chi_2(x)f_u(x)\in C^4(\mathcal{U})$. It has been proved in Bai, Wang and Zhou \cite{BWZ} that for $C^4(\mathcal{U})$ function supported on $\mathcal{U}_\epsilon$, one has
\begin{eqnarray}
\mathbb{E}\mathcal{L}_n[\chi_2(x)f_u(x)]=n\int \chi_2(x)f_u(x)dF_{sc}(x)+o(1). \label{2.9}
\end{eqnarray}
For $\mathbb{E}\mathcal{L}_n[\chi_1(x)f_u(x)]$, we use the following asymptotic formula proved in Ercolani and McLaughlin \cite{EM},
\begin{eqnarray*}
\mathcal{K}_n(x,x)=\frac{n}{2\pi}\sqrt{4-x^2}+\frac{1}{4\pi}(\frac{1}{x-2}-\frac{1}{x+2})\cos [\frac{n}{2\pi}\int_x^2\sqrt{4-y^2}dy]+\mathcal{O}(n^{-1})
\end{eqnarray*}
for $x\in[-2+\delta,2-\delta]$ with any fixed $\delta>0$. Thus by using (\ref{2.25}) one has
\begin{eqnarray*}
\mathbb{E}\mathcal{L}_n[\chi_1(x)f_u(x)]=n\int \chi_1(x) f_u(x)dF_{sc}(x)+\varepsilon_n,
\end{eqnarray*}
where
\begin{eqnarray*}
\varepsilon_n=\int_{-\frac 32+\frac u4}^u \chi_1(x) f_u(x)\frac{1}{4\pi}(\frac{1}{x-2}-\frac{1}{x+2})\cos [\frac{n}{2\pi}\int_x^2\sqrt{4-y^2}dy]dx+\mathcal{O}(n^{-1}).
\end{eqnarray*}
By integration by parts, we can easily get that
\begin{eqnarray}
\varepsilon_n=\mathcal{O}(n^{-1}). \label{2.10}
\end{eqnarray}
Thus combining (\ref{2.9}) with (\ref{2.10}) we can show (\ref{2.7}).

Now we prove (iii) of Theorem \ref{th.1}. By definition, $\gamma_{k/n}$ is the $k$-th $n$-quantile of the semicircle law, i.e.
\begin{eqnarray*}
\frac{1}{2\pi}\int_{-2}^{\gamma_{k/n}}\sqrt{4-x^2}dx=\frac{k}{n}.
\end{eqnarray*}
We decompose $\mathcal{B}_n[f;k]$ as
\begin{eqnarray*}
\mathcal{B}_n[f;k]=\sum_{l=1}^k(f(\lambda_l)-f(\gamma_{k/n}))+kf(\gamma_{k/n}).
\end{eqnarray*}
In order to avoid working on the ordered eigenvalues, we introduce a proxy of $\mathcal{B}_n[f;k]$ as
\begin{eqnarray*}
\widehat{\mathcal{B}}_n[f;k]
&=&\sum_{l=1}^n(f(\lambda_l)-f(\gamma_{k/n}))\mathbf{1}_{\{\lambda_l\leq \gamma_{k/n}\}}+kf(\gamma_{k/n})\nonumber\\
&=&\sum_{l=1}^{nF_n(\gamma_{k/n})}(f(\lambda_l)-f(\gamma_{k/n}))+kf(\gamma_{k/n}).
\end{eqnarray*}
Let $a=|nF_n(\gamma_{k/n})-k|$. By using the rigidity property in Lemma \ref{lem.9}, one has with overwhelming probability
$$a=|nF_n(\gamma_{k/n})-nF_{sc}(\gamma_{k/n})|\leq(\log n)^{C\log\log n}.$$
Furthermore, we also have
\begin{eqnarray}
|\mathcal{B}_n[f;k]-\widehat{\mathcal{B}}_n[f;k]|\leq a\max_{l\in\{k-a,k+a\}} |f(\lambda_l)-f(\gamma_{k/n})|=\mathcal{O}(\frac{(\log n)^{C\log\log n}}{n}) \label{2.18}
\end{eqnarray}
with overwhelming probability.  Therefore, we only have to prove the central limit theorem for $\widehat{\mathcal{B}}_n[f;k]$.

Observe that
\begin{eqnarray}
\widehat{\mathcal{B}}_n[f;k]=\sum_{l=1}^nf_{\gamma_{k/n}}(\lambda_l)+kf(\gamma_{k/n})= \mathcal{L}_n[f_{\gamma_{k/n}}]+kf(\gamma_{k/n}). \label{2.19}
\end{eqnarray}
Moreover, since we assume $f$ is compactly supported on $\mathcal{U}_\epsilon$, we have
\begin{eqnarray*}
\mathbb{E}\mathcal{L}_n[f_{\gamma_{k/n}}]=m[f;\gamma_{k/n}]+o(1);\quad \mathbb{E}\mathcal{L}_n[f_{\gamma_{y}}]=m[f;\gamma_{y}]+o(1)
\end{eqnarray*}
as shown in (\ref{2.7}).
Thus to prove (iii), it suffices to show as $n\rightarrow\infty$
\begin{eqnarray}
|\mathbb{E}e^{-ix\mathcal{L}^\circ_n[f_{\gamma_{k/n}}]}-\mathbb{E}e^{-ix\mathcal{L}^\circ_n[f_{\gamma_{y}}]}|\rightarrow 0 \label{2.21}
\end{eqnarray}
for any fixed $x$. To see (\ref{2.21}), we note that
\begin{eqnarray*}
&&|\mathbb{E}e^{-ix\mathcal{L}^\circ_n[f_{\gamma_{k/n}}]}-\mathbb{E}e^{-ix\mathcal{L}^\circ_n[f_{\gamma_{y}}]}|\nonumber\\
&&\leq |x|\mathbb{E}|\mathcal{L}^\circ_n[f_{\gamma_{k/n}}]-\mathcal{L}^\circ_n[f_{\gamma_{y}}]|\nonumber\\
&&\leq|x|Var^{1/2}\{\mathcal{L}_n[f_{\gamma_{k/n}}-f_{\gamma_y}]\}.
\end{eqnarray*}
Thus it remains to verify
\begin{eqnarray}
Var\mathcal{L}_n[f_{\gamma_{k/n}}-f_{\gamma_y}]\rightarrow 0. \label{2.22}
\end{eqnarray}
 To show this, we will rely on the following lemma whose proof will be postponed to the end of this section. It will be clear that the following lemma is also crucial to our proof of tightness for $\mathcal{S}^\circ_n[f;t]$ in the next section.
\begin{lem}\label{lem.7}Suppose that $\varphi$ is a Lipschitz function on $\mathbb{R}$ with Lipschitz constant $L$. Moreover, we assume that there exists an interval $I=[a,b]\subset[-2+\delta,2-\delta]$ ($a,b$ may be $n$-dependent) such that
\begin{eqnarray*}
\varphi(\lambda)=\varphi(a), \quad \lambda\leq a;\quad \varphi(\lambda)=\varphi(b),\quad \lambda\geq b.
\end{eqnarray*}
Then for GUE, we have
\begin{eqnarray*}
Var\mathcal{L}_n[\varphi]\leq C(b-a)^2(|\log(b-a)|+1)+\mathcal{O}(e^{-cn})
\end{eqnarray*}
with some positive constants $C=:C(L)$ and $c$ independent of $a$ and $b$.
\end{lem}
Now we proceed to the proof of (iii) of Theorem \ref{th.1}. For convenience, we assume $k/n\geq y$. The opposite case is just analogous. Let $\varphi=f_{\gamma_{k/n}}-f_{\gamma_y}$. By definition, we note that $\varphi(\lambda)$ equals to $f(\gamma_y)-f(\gamma_{k/n})$ for $\lambda\leq\gamma_y $ and $0$ for $\lambda\geq \gamma_{k/n}$. Thus by Lemma \ref{lem.7} and the assumption that $k/n\rightarrow y$, we have (\ref{2.22}). Thus (\ref{2.21}) holds. So we conclude the proof of (iii) by using (\ref{2.21}) and (ii) of Theorem \ref{th.1}.
\end{proof}

\begin{proof}[Proof of Lemma \ref{lem.7}]
From the determinantal structure of the j.p.d. of eigenvalues for GUE, with (\ref{2.26}) one has
\begin{eqnarray*}
Var\mathcal{L}_n[\varphi]=\int_\mathbb{R}\int_\mathbb{R}|\frac{\varphi(\lambda)-\varphi(\mu)}{\lambda-\mu}|^2\mathcal{V}_n(\lambda,\mu)d\lambda d\mu,
\end{eqnarray*}
where
\begin{eqnarray*}
\mathcal{V}_n(\lambda,\mu)&=&\frac12\mathcal{K}^2_n(\lambda,\mu)(\lambda-\mu)^2\nonumber\\
&=&\frac12(\sqrt{n}\psi_n(\sqrt{n}x)\psi_{n-1}(\sqrt{n}y)-\sqrt{n}\psi_{n-1}(\sqrt{n}x)\psi_n(\sqrt{n}y))^2.
\end{eqnarray*}
By assumption, we can split the integral into four parts
\begin{eqnarray*}
Var\mathcal{L}_n[\varphi]&=& V_1+V_2+V_3+V_4
\end{eqnarray*}
where
\begin{eqnarray*}
&&V_1=:\int_{a}^{b}\int_{a}^{b}|\frac{\varphi(\lambda)-\varphi(\mu)}{\lambda-\mu}|^2\mathcal{V}_n(\lambda,\mu)d\lambda d\mu\nonumber\\
&&V_2=:2\int_{a}^{b}\int_{-\infty}^{a}|\frac{\varphi(\lambda)-\varphi(\mu)}{\lambda-\mu}|^2\mathcal{V}_n(\lambda,\mu)d\lambda d\mu\nonumber\\
&&V_3=:2\int_{a}^{b}\int_{b}^{\infty}|\frac{\varphi(\lambda)-\varphi(\mu)}{\lambda-\mu}|^2\mathcal{V}_n(\lambda,\mu)d\lambda d\mu\nonumber\\
&&V_4=:2\int_{-\infty}^{a}\int_{b}^{\infty}|\frac{\varphi(\lambda)-\varphi(\mu)}{\lambda-\mu}|^2\mathcal{V}_n(\lambda,\mu)d\lambda d\mu.
\end{eqnarray*}
Note that $\varphi$ is Lipschitz. And it is well known that $\mathcal{V}_n(x,y)$ is bounded in $\mathbb{R}^2$. Moreover, by (\ref{2.29}) one sees that $\mathcal{V}_n(x,y)$ is exponentialy decaying in $x,y$ as $|x|$ or $|y|$ larger than $2+\delta$. Thus we immediately get that
\begin{eqnarray*}
&&V_1\leq C(b-a)^2,\nonumber\\
&&V_2=:2\int_{a}^{b}\int_{-2-\delta}^{a}|\frac{\varphi(\lambda)-\varphi(\mu)}{\lambda-\mu}|^2\mathcal{V}_n(\lambda,\mu)d\lambda d\mu+\mathcal{O}(e^{-cn}),\nonumber\\
&&V_3=:2\int_{a}^{b}\int_{b}^{2+\delta}|\frac{\varphi(\lambda)-\varphi(\mu)}{\lambda-\mu}|^2\mathcal{V}_n(\lambda,\mu)d\lambda d\mu+\mathcal{O}(e^{-cn}),\nonumber\\
&&V_4=:2\int_{-2-\delta}^{a}\int_{b}^{2+\delta}|\frac{\varphi(\lambda)-\varphi(\mu)}{\lambda-\mu}|^2\mathcal{V}_n(\lambda,\mu)d\lambda d\mu+\mathcal{O}(e^{-cn}).
\end{eqnarray*}

Now we estimate $V_2$. Note that
\begin{eqnarray*}
&&V_2=2\int_{a}^{b}\int_{2a-b}^{a}|\frac{\varphi(\lambda)-\varphi(\mu)}{\lambda-\mu}|^2\mathcal{V}_n(\lambda,\mu)d\lambda d\mu\nonumber\\
&&+2\int_{a}^{b}\int_{-2-\delta}^{2a-b}|\frac{\varphi(\lambda)-\varphi(\mu)}{\lambda-\mu}|^2\mathcal{V}_n(\lambda,\mu)d\lambda d\mu+\mathcal{O}(e^{-cn}).\nonumber\\
\end{eqnarray*}
Observe that the first term on the right hand side of the above equality can be bounded by
$C(b-a)^2,$
and the second term can be bounded as
\begin{eqnarray*}
&&\int_{a}^{b}\int_{-2-\delta}^{2a-b}|\frac{\varphi(\lambda)-\varphi(\mu)}{\lambda-\mu}|^2\mathcal{V}_n(\lambda,\mu)d\lambda d\mu \nonumber\\
&&\leq C(b-a)^2\int_{a}^{b}\int_{-2-\delta}^{2a-b}|\frac{1}{\lambda-\mu}|^2d\lambda d\mu \nonumber\\
&&\leq C(b-a)^2(|\log(b-a)|+1).
\end{eqnarray*}
Thus we have
\begin{eqnarray*}
V_2\leq C(b-a)^2(|\log(b-a)|+1)+\mathcal{O}(e^{-cn}).
\end{eqnarray*}
Analogously, one can also get that
\begin{eqnarray*}
V_3,V_4\leq C(b-a)^2(|\log(b-a)|+1)+\mathcal{O}(e^{-cn}).
\end{eqnarray*}
Thus we conclude the proof.
\end{proof}
\section{Partial sum process for GUE}
In this section, we provide the proof of Theorem \ref{th.3}. Thus we have to verify the finite dimensional convergence and the tightness of the sequence $\{\mathcal{S}^\circ_n[f;t]; t\in [\delta,1-\delta]\}$.

At first, we extend the discussion in the last section to show that the finite dimensional convergence of the process $\{\mathcal{S}^\circ_n[f;t]; t\in [\delta,1-\delta]\}$. We formulate the result as the following lemma.
\begin{lem} \label{lem.8}Under the assumptions of Theorem \ref{th.3}, for any fixed positive integer $r$ and points $t_1,\cdots,t_r\in[\delta,1-\delta]$, and for any fixed numbers $\alpha_1,\cdots,\alpha_r\in \mathbb{R}$, we have
\begin{eqnarray*}
\sum_{l=1}^r\alpha_l\mathcal{S}^\circ_n[f;t_l]\stackrel{d}\longrightarrow \sum_{l=1}^r\alpha_l\mathcal{S}[f;t_l].
\end{eqnarray*}
\end{lem}
\begin{proof} Below we set $k_l=\lfloor nt_l\rfloor$ and $u(t_l)=k_l/n$. At first, we claim that for integer $l\in [\delta n,(1-\delta)n]$,
\begin{eqnarray}
Var\{f({\lambda_l})\}\leq C\frac{\log n}{n^2} \label{4.1}
\end{eqnarray}
for $f$ obeying the assumptions in Theorem \ref{th.3}. Here the constant $C$ depends only on $\delta$ and the function $f$. To show (\ref{4.1}), we recall the bump function $\chi_\epsilon$ and the corresponding truncated function $\bar{f}_\epsilon(x)$ defined in Section 3. Then we put
\begin{eqnarray*}
\tilde{f}_\epsilon(x)=f(x)-\bar{f}_\epsilon(x),
\end{eqnarray*}
which vanishes when $|x|\leq 2+\epsilon/2$. Therefore, one has
\begin{eqnarray*}
Var\{f({\lambda_l})\}\leq 2Var\{\bar{f}_\epsilon({\lambda_l})\}+2\mathbb{E}(\tilde{f}_\epsilon(\lambda_l))^2.
\end{eqnarray*}
It follows from Lemma \ref{lem.10} that
\begin{eqnarray*}
Var\{\bar{f}_\epsilon({\lambda_l})\}\leq \sup_{x}|\bar{f}'_\epsilon(x)|^2Var\{\lambda_l\}\leq C\frac{\log n}{n^2}.
\end{eqnarray*}
Besides, by (\ref{2.33}) and the assumptions on $f(x)$, we also have
\begin{eqnarray*}
\mathbb{E}(\tilde{f}_\epsilon(\lambda_l))^2=\mathcal{O}(e^{-cn})\leq C\frac{\log n}{n^2}
\end{eqnarray*}
for sufficiently large $n$. Thus we have (\ref{4.1}).
Consequently, we have
\begin{eqnarray*}
\sum_{l=1}^r\alpha_lf^\circ(\lambda_{k_l+1})\stackrel{P}\longrightarrow 0.
\end{eqnarray*}
Then it remains to show
\begin{eqnarray*}
\sum_{l=1}^r\alpha_lB^\circ_n[f; k_l]\stackrel{d}\longrightarrow \sum_{l=1}^r\alpha_l\mathcal{S}[f;t_l].
\end{eqnarray*}
Using (\ref{2.18}) and (\ref{2.19}), it suffices to prove
\begin{eqnarray*}
\sum_{l=1}^r\alpha_l\mathcal{L}^\circ_n[f_{\gamma_{u(t_l)}}]\stackrel{d}\longrightarrow \sum_{l=1}^r\alpha_l\mathcal{S}[f;t_l].
\end{eqnarray*}
By the fact that $u(t_l)\rightarrow t_l$ and an routine discussion as that for (\ref{2.21}), we can reduce the problem to show that
\begin{eqnarray*}
\sum_{l=1}^r\alpha_l\mathcal{L}^\circ_n[f_{\gamma_{t_l}}]\stackrel{d}\longrightarrow \sum_{l=1}^r\alpha_l\mathcal{S}[f;t_l].
\end{eqnarray*}
Note that
\begin{eqnarray*}
\sum_{l=1}^r\alpha_l\mathcal{L}^\circ_n[f_{\gamma_{t_l}}]=\mathcal{L}_n^\circ[\sum_{l=1}^r\alpha_l f_{\gamma_{t_l}}].
\end{eqnarray*}
Observe that $\sum_{l=1}^r\alpha_l f_{\gamma_{t_l}}$ is a continuous function with $r$ possibly non differentiable points $t_1, \cdots, t_r$. Now we choose $r$ interval $I_1,\cdots,I_r$ containing $\gamma_{t_1},\cdots,\gamma_{t_r}$ respectively with lengths $|I_l|\leq n^{-1/2-c}$ for some small positive number $c$ and all $l=1,\cdots,r$. We define a smooth modification function $g_{t_1,\cdots,t_r}$ which coincides with $\sum_{l=1}^r\alpha_l f_{\gamma_{t_l}}$ on $\mathbb{R}\setminus\cup_{l=1}^rI_l$ and obeys the condition
\begin{eqnarray*}
|\frac{d^{k+1}}{dx^{k+1}}g_{t_1,\cdots,t_r}(x)|=\mathcal{O}(n^{k(1/2+c)}),k=0,1,2.
\end{eqnarray*}
By a similar relation to (\ref{2.8}), we only have to prove
\begin{eqnarray}
\mathcal{L}_n^\circ[g_{t_1,\cdots,t_r}]\stackrel{d}\longrightarrow \sum_{l=1}^r\alpha_l\mathcal{S}[f;t_l]. \label{2.23}
\end{eqnarray}
The proof of (\ref{2.23}) is easy to carry out by using Lytova and Pastur's method again as that in Section 3. Thus we can finally get that $\sum_{l=1}^r\alpha_l\mathcal{S}^\circ_n[f;t_l]$ weakly converges to a Gaussian distribution with mean zero and variance $V_{GUE}[\sum_{l=1}^r\alpha_l f_{\gamma_{t_l}}]$. Consequently, we can conclude the proof of Lemma \ref{lem.8}.
\end{proof}

It remains to show that the sequence $\{\mathcal{S}^\circ_n[f;t];t\in[\delta,1-\delta]\}$ is tight. To this end, we will use Theorem 12.3 (p. 95) of Billingsley \cite{Bi}. According to this theorem, we need to verify

\emph{({\bf{i}}): Tightness at any point in $[\delta,1-\delta]$.}

\emph{({\bf{ii}}): For arbitrary $s,t\in[\delta,1-\delta]$ and $n$ sufficiently large
\begin{eqnarray*}
\mathbb{E}|\mathcal{S}^\circ_n[f;t]-\mathcal{S}^\circ_n[f;s]|^2\leq C|t-s|^\alpha
\end{eqnarray*}
for some constant $C>0$,and $\alpha>1$ which are independent of $t,s$.}

Note that ({\bf{i}}) is obvious. Thus it suffices to show ({\bf{ii}}). Set $\eta_n=\log^{-L}n$ for some constant $L$ large enough. Without loss of generality, we always assume that $s\leq t$ below. We separate the issue into three cases: $0\leq t-s\leq n^{-1}$, $n^{-1}\leq t-s\leq \eta_n$ and $t-s\geq \eta_n$.

For $t-s\leq n^{-1}$, one has $\lfloor nt\rfloor=\lfloor ns\rfloor$ or $\lfloor nt\rfloor=\lfloor ns\rfloor+1$. When $\lfloor nt\rfloor=\lfloor ns\rfloor$, we have
\begin{eqnarray*}
\mathbb{E}|\mathcal{S}^\circ_n[f;t]-\mathcal{S}^\circ_n[f;s]|^2=\mathbb{E}|n(t-s)f^\circ(\lambda_{\lfloor ns\rfloor+1})|^2\leq C|t-s|^2\log n.
\end{eqnarray*}
In the last step above we used the estimation (\ref{4.1}). And the positive constant $C$ only depends on $\delta$ and the test function $f$ (Such a dependence will not be mentioned repeatedly below for simplicity). When $\lfloor nt\rfloor=\lfloor ns\rfloor+1$, one has
\begin{eqnarray*}
\mathbb{E}|\mathcal{S}^\circ_n[f;t]-\mathcal{S}^\circ_n[f;s]|^2=\mathbb{E}|(nt-\lfloor nt\rfloor)f^\circ(\lambda_{\lfloor nt\rfloor+1})+(\lfloor nt\rfloor-ns)f^\circ(\lambda_{\lfloor nt\rfloor})|^2.
\end{eqnarray*}
Note that when $\lfloor nt\rfloor=\lfloor ns\rfloor+1$, one has
\begin{eqnarray*}
0\leq nt-\lfloor nt\rfloor, \lfloor nt\rfloor-ns\leq nt-ns.
\end{eqnarray*}
Consequently, we have
\begin{eqnarray*}
&&\mathbb{E}|\mathcal{S}^\circ_n[f;t]-\mathcal{S}^\circ_n[f;s]|^2\nonumber\\
&&\leq Cn^2(t-s)^2(Var\{f(\lambda_{\lfloor nt\rfloor+1})\}
+Var\{f(\lambda_{\lfloor nt\rfloor})\})\nonumber\\
&&\leq C|t-s|^2\log n .
\end{eqnarray*}
Since $0\leq t-s\leq n^{-1}$, for $n$ large enough, we always have
\begin{eqnarray*}
|t-s|^2\log n\leq |t-s|^{3/2}.
\end{eqnarray*}

For $n^{-1}\leq t-s\leq \eta_n$, one has
\begin{eqnarray*}
\mathbb{E}|\mathcal{S}^\circ_n[f;t]-\mathcal{S}^\circ_n[f;s]|^2
&\leq& 2\mathbb{E}|\mathcal{B}^\circ_n[f,\lfloor nt\rfloor]-\mathcal{B}^\circ_n[f,\lfloor ns\rfloor]|^2\nonumber\\
&&+2\mathbb{E}|(nt-\lfloor nt\rfloor)f^\circ(\lambda_{\lfloor nt\rfloor+1})-(ns-\lfloor ns\rfloor)f^\circ(\lambda_{\lfloor ns\rfloor+1})|^2\nonumber\\
&\leq&2(\lfloor nt\rfloor-\lfloor ns\rfloor)\sum_{i=\lfloor ns\rfloor}^{\lfloor nt\rfloor}Var\{f(\lambda_i)\}\nonumber\\
&&+2\mathbb{E}|(nt-\lfloor nt\rfloor)f^\circ(\lambda_{\lfloor nt\rfloor+1})-(ns-\lfloor ns\rfloor)f^\circ(\lambda_{\lfloor ns\rfloor+1})|^2\nonumber\\
&\leq& C (t-s)^2\log n.
\end{eqnarray*}
In the above second inequality we have used the basic relation
\begin{eqnarray*}
Var\{\xi_1+\cdots+\xi_m\}\leq m\sum_{l=1}^m Var\{\xi_l\}.
\end{eqnarray*}
Clearly, when the constant $L$ in the definition of $\eta_n$ is chosen to be large enough, for $n^{-1}\leq t-s\leq \eta_n$ we have
\begin{eqnarray*}
(t-s)^2\log n\leq (t-s)^{3/2}.
\end{eqnarray*}

For the last case $t-s\geq \eta_n$, it suffices to show that
\begin{eqnarray*}
\mathbb{E}|\mathcal{B}^\circ_n[f,\lfloor nt\rfloor]-\mathcal{B}^\circ_n[f,\lfloor ns\rfloor]|^2\leq C|t-s|^{3/2}.
\end{eqnarray*}
By (\ref{2.18}) and (\ref{2.19}) we need to prove
\begin{eqnarray}
Var|\mathcal{L}_n[f_{\gamma_{u(t)}}-f_{\gamma_{u(s)}}]|^2\leq C|t-s|^{3/2}, \label{2.24}
\end{eqnarray}
where $u(t)=\lfloor nt\rfloor/n$. Note that (\ref{2.24}) follows from Lemma \ref{lem.7} immediately. So Theorem \ref{th.3} follows.
\section{CLTs for Wigner matrices}
As shown in Section 3, Lytova and Pastur's original proof in \cite{LP} for the GUE case can be easily modified to adapt to our case. However, for more general complex Wigner matrix, higher order derivatives of $g_u(x)$ will be involved if we proceed to pursue the discussion in \cite{LP} (see (3.49) of \cite{LP} for instance). But those derivatives in (\ref{2.17}) will not be small enough for the strategy in \cite{LP}. Moreover, the results in \cite{LP} do not provide the asymptotic estimation of the expectation.

Motivated by the recent articles \cite{TV2} and \cite{TV3}, we will establish a comparison theorem for the linear eigenvalue statistics
with a certain class of test functions in this section. As an application, we use the comparison theorem to extend Theorem \ref{th.1} to general complex Wigner matrices case. Moreover, such a comparison theorem will also be used in the next section to prove Theorem \ref{th.4}.

At first, we define the set of $n$-dependent real functions $\mathcal{F}_n^m$ for some fixed positive integer $m$. We say a function $\varphi\in \mathcal{F}_n^m$ if and only if $\varphi=:\varphi_n$ satisfies the following assumptions ({\bf{a}}) and ({\bf{b}}).

({\bf{a}}): \emph{$\varphi\in C^4(\mathcal{U})$ is compactly supported on $\mathcal{U}_{\epsilon}$ and
\begin{eqnarray*}
|\varphi^{(\alpha)}(x)|\leq C, \quad \alpha=0,1
\end{eqnarray*}
 with some positive constant $C$ independent of $n$.}

({\bf{b}}): \emph{There exist $m$ intervals $I_1,\cdots, I_m\in[-2+\delta,2-\delta]$ with length $|I_l|\leq n^{-1-c_1}$ for some small $c_1>0$ and all $l=1,\cdots,m$, such that for $x\in\mathcal{U}\setminus {\cup_{l=1}^m I_l}$ one has
\begin{eqnarray*}
|\varphi^{(\alpha)}(x)|\leq C,\quad \alpha=2,3,4
\end{eqnarray*}
 with some positive constant $C$ independent of $n$}.

Our main tool to extend the CLTs from GUE to general Wigner matrices is the following comparison theorem for linear eigenvalue statistics.
\begin{thm} \label{th.3.1} Let $M_n=\frac{1}{\sqrt{n}}(w_{jk})_{j,k=1}^n$ and $M_n'=\frac{1}{\sqrt{n}}(w'_{jk})_{j,k=1}^n$ be two Wigner matrices satisfying Condition $\mathbf{C_0}$. We assume $M_n$ and $M'_n$ match to the fourth order off the diagonal and to the second order on the diagonal. Moreover, the magnitudes of $w_{jk}$ and $w'_{jk}$ ($1\leq j,k\leq n$) are uniformly bounded by $n^{\mathcal{O}(c_0)}$ for some sufficiently small but fixed $c_0>0$. Let
$G:\mathbb{R}\rightarrow\mathbb{R}$ obey the derivative bounds
\begin{eqnarray}
|\frac{d^jG(x)}{dx^j}|=\mathcal{O}(n^{c_0}) \label{5.3}
\end{eqnarray}
for $0\leq j\leq 5$. If $\varphi\in \mathcal{F}_n^m$ for some fixed integer $m$, then we have
\begin{eqnarray*}
\mathbb{E}G(\mathcal{L}^M_n[\varphi])-\mathbb{E}G(\mathcal{L}^{M'}_n[\varphi])=\mathcal{O}(n^{-c})
\end{eqnarray*}
for some fixed $c>0$. Here $\mathcal{L}^M_n[\varphi]$ (resp. $\mathcal{L}^{M'}_n[\varphi]$) represents the linear eigenvalue statistic of $M_n$ (resp. $M_n'$) with the test function $\varphi$.
\end{thm}
Using the terminology in \cite{TV2}, we say a statistic $S(M_n)$ that can depend on $M_n$ or $M'_n$ \emph{highly insensitive} if one has
\begin{eqnarray*}
|S(M_n)-S(M'_n)|=\mathcal{O}(n^{-c}).
\end{eqnarray*}
for some fixed $c>0$. Thus Theorem \ref{th.3.1} asserts that $\mathbb{E}G(\mathcal{L}_n^M[\varphi])$ is highly insensitive for $\varphi\in\mathcal{F}_n^m$.
To show this, our strategy is to represent the linear eigenvalue statistics by the Stieltjes transform of the ESD defined in Section 2. Then the Lindeberg swapping argument for the Stieltjes transform which was well developed in recent work such as \cite{EYY} and \cite{TV2} can be applied. To this end, we use the following
Helffer-Sj\"{o}strand formula.
\begin{lem}[Helffer-Sj\"{o}strand formula] Suppose that $\varphi:\mathbb{R}\rightarrow\mathbb{R}$ be a $C^{k+1}(\mathbb{R})$ function with a compact support. Let $\sigma(y)\in C^{\infty}(\mathbb{R})$ be a cut off function such that $\sigma(y)=1$ for $|y|\leq1/2$ and $\sigma(y)=0$ for $|y|\geq 1$ with bounded derivatives.
Define the smooth extension $\tilde{\varphi}:\mathbb{C}\rightarrow\mathbb{C}$ of $\varphi$ by
\begin{eqnarray*}
\tilde{\varphi}(z):=\left(\sum_{j=0}^{k}\frac{\varphi^{(j)}(x)(iy)^j}{j!}\right)\sigma(y),
\end{eqnarray*}
where $z=x+iy$. Then for any self-adjoint operator $X$, one has
\begin{eqnarray*}
\varphi(X)=\frac1\pi\int_{\mathbb{R}^2}\frac{\partial \tilde{\varphi}}{\partial \bar{z}}\frac{1}{X-z}dxdy,
\end{eqnarray*}
where
\begin{eqnarray*}\frac{\partial\tilde{\varphi}}{\partial\bar{z}}
:=\frac12\left(\frac{\partial\tilde{\varphi}}{\partial x}+i\frac{\partial\tilde{\varphi}}{\partial y}\right).
\end{eqnarray*}
\end{lem}
\begin{rem}We refer to Davies' book \cite{D} for more details on the Helffer-Sj\"{o}strand formula. Moreover, in the literature of RMT,
one can also see \cite{ESY2} and \cite{RRS} for references.
\end{rem}
By definition, one can calculate
\begin{eqnarray*}
\frac{\partial\tilde{\varphi}}{\partial\bar z}=\frac12\left(\sum_{j=0}^{k}\frac{\varphi^{(j)}(x)(iy)^j}{j!}\right)i\sigma'(y)
+\frac{1}{2k!}\varphi^{(k+1)}(x)(iy)^k\sigma(y).
\end{eqnarray*}

Below we will restrict to the case of $m=1$ for ease of presentation.  It will be clear that the proof can be extended straightforward to the case of $m>1$ but fixed. For simplicity, we will denote $I_1$ and $\mathcal{F}_n^1$ by $I$ and $\mathcal{F}_n$ respectively. For $\varphi\in\mathcal{F}_n$, we denote $I=:[a,b]$ with $a=a_n,b=b_n$ such that $|b-a|\leq n^{-1-c_1}$. By the Helffer-Sj\"{o}strand formula we can write
\begin{eqnarray}
\mathcal{L}_n[\varphi]&=&\frac 1\pi\sum_{l=1}^n\int_{\mathbb{R}^2}\frac{\partial \tilde{\varphi}}{\partial \bar{z}}\frac{1}{\lambda_l-z}dxdy\nonumber\\
&=&\frac{1}{2\pi}\sum_{l=1}^n\int_{\mathbb{R}^2}\left(\sum_{j=0}^{3}\frac{\varphi^{(j)}(x)(iy)^j}{j!}\right)i\sigma'(y)\frac{1}{\lambda_l-z}dxdy\nonumber\\
&&~~~+\frac{1}{2\pi}\sum_{l=1}^n\int_{\mathbb{R}^2}\frac{1}{3!}\varphi^{(4)}(x)(iy)^3\sigma(y)\frac{1}{\lambda_l-z}dxdy\nonumber\\
&=&\frac{1}{2\pi}\sum_{l=1}^n\int_{-2-\epsilon}^{2+\epsilon}\left(\int_{-1}^{-1/2}+\int_{1/2}^{1}\right)\left(\sum_{j=0}^3\frac{\varphi^{(j)}(x)(iy)^j}{j!}\right)i\sigma'(y)
\frac{1}{\lambda_l-z}dxdy\nonumber\\
&&+\frac{1}{2\pi}\sum_{l=1}^n\left(\int_{b}^{2+\epsilon}+\int_{-2-\epsilon}^{a}\right)\int_{-1}^1\frac{1}{3!}\varphi^{(4)}(x)(iy)^3\sigma(y)\frac{1}{\lambda_l-z}dxdy \nonumber\\
&&+\frac{1}{2\pi}\sum_{l=1}^n\int_{a}^{b}\int_{-1}^1\frac{1}{3!}\varphi^{(4)}(x)(iy)^3\sigma(y)\frac{1}{\lambda_l-z}dxdy\nonumber\\
&:=&A_1+A_2+A_3. \label{1.27}
\end{eqnarray}

Observe that  $|y|\geq 1/2$ in the integral region of the first term $A_1$. Because $(\lambda_l-z)^{-1}$ is analytic in this region, we can use integration by parts. It is not difficult to derive that
\begin{eqnarray*}
A_1=\frac{n}{2\pi}\int_{-2-\epsilon}^{2+\epsilon}\left(\int_{-1}^{-1/2}+\int_{1/2}^{1}\right)\left(\sum_{j=0}^{3}\frac{1}{j!}i[y^j\sigma'(y)]^{(j)}\right)
\varphi(x)s_n(x+iy)dxdy.
\end{eqnarray*}

For the term $A_2$, we note that since $x\in \mathcal{U}_\epsilon\setminus I$, $|\varphi^{(4)}(x)|\leq C$ with some positive constant $C$ independent of $n$ by assumption. Moreover, we always have the elementary inequality
\begin{eqnarray}
|(\lambda_l-z)^{-1}|\leq y^{-1}. \label{5.4}
\end{eqnarray}
Let $y_0=n^{-1-c_0}$, we decompose the integral region in $A_2$ into two parts: $|y|\leq y_0$ and $|y|>y_0$. Then (\ref{5.4}) implies that
\begin{eqnarray*}
\left|\frac{1}{2\pi}\sum_{l=1}^n\left(\int_{b}^{2+\epsilon}+\int_{-2-\epsilon}^{a}\right)\int_{|y|\leq y_0}\frac{1}{3!}\varphi^{(4)}(x)(iy)^3\sigma(y)\frac{1}{\lambda_l-z}dxdy\right|\leq Cn^{-2-3c_0}.
\end{eqnarray*}
Therefore, we can write
\begin{eqnarray*}
A_2=\frac{n}{2\pi}\left(\int_{b}^{2+\epsilon}+\int_{-2-\epsilon}^{a}\right)\int_{|y|> y_0}\frac{1}{3!}\varphi^{(4)}(x)(iy)^3\sigma(y)s_n(z)dxdy+\mathcal{O}(n^{-2-3c_0}).
\end{eqnarray*}

For the third term $A_3$, we will condition on the event $N_n(I)=0$. It is clear that if there is no eigenvalue in the interval $I=[a,b]$, then $(\lambda_l-z)^{-1}$ is continuously differentiable w.r.t $x$ and $y$ in the integral region $I\times[-1,1]$. Consequently, when $N_n(I)=0$, we can apply integration by parts to the term $A_3$ and obtain
\begin{eqnarray*}
&&A_3=\frac{1}{2\pi}\frac{1}{3!}\sum_{l=1}^n\sum_{\alpha=0}^2i^{3-\alpha}\int_{-1}^1(y^3\sigma(y))^{(\alpha)}
\left(\frac{\varphi^{(3-\alpha)}(b)}{\lambda_l-b-iy}
-\frac{\varphi^{(3-\alpha)}(a)}{\lambda_l-a-iy}\right)dy\nonumber\\
&&+\frac{1}{2\pi}\frac{1}{3!}\sum_{l=1}^n\int_a^b\int_{-1}^1\varphi'(x)(y^3\sigma(y))^{(3)}\frac{1}{\lambda_l-z}dxdy\nonumber\\
&&=\frac{n}{2\pi}\frac{1}{3!}\sum_{\alpha=0}^2i^{3-\alpha}\int_{|y|> n^{-4}}(y^3\sigma(y))^{(\alpha)}\left(\varphi^{(3-\alpha)}(b)s_n(b+iy)-\varphi^{(3-\alpha)}(a)s_n(a+iy)\right)dy\nonumber\\
&&+\frac{1}{2\pi}\frac{1}{3!}\sum_{l=1}^n\int_a^b\int_{-1}^1\varphi'(x)(y^3\sigma(y))^{(3)}\frac{1}{\lambda_l-z}dxdy+\mathcal{O}(n^{-3}).\nonumber\\
\end{eqnarray*}
In the last step, we used (\ref{5.4}) again to assert
\begin{eqnarray*}
\frac{1}{2\pi}\frac{1}{3!}\sum_{l=1}^n\sum_{\alpha=0}^2i^{3-\alpha}\int_{|y|\leq n^{-4}}(y^3\sigma(y))^{(\alpha)}
\left(\frac{\varphi^{(3-\alpha)}(b)}{\lambda_l-b-iy}
-\frac{\varphi^{(3-\alpha)}(a)}{\lambda_l-a-iy}\right)dy=\mathcal{O}(n^{-3}).
\end{eqnarray*}
Consequently, when $N_n(I)=0$, we obtain
\begin{eqnarray}
&&\mathcal{L}_n[\varphi]=\frac{n}{2\pi}\int_{-2-\epsilon}^{2+\epsilon}\left(\int_{-1}^{-1/2}
+\int_{1/2}^{1}\right)\left(\sum_{j=0}^{3}\frac{1}{j!}i(y^j\sigma'(y))^{(j)}\right)\varphi(x)s_n(x+iy)dxdy\nonumber\\
&&+\frac{n}{2\pi}\left(\int_{b}^{2+\epsilon}+\int_{-2-\epsilon}^{a}\right)\int_{|y|> y_0}\frac{1}{3!}\varphi^{(4)}(x)(iy)^3\sigma(y)s_n(x+iy)dxdy\nonumber\\
&&+\frac{n}{2\pi}\frac{1}{3!}\sum_{\alpha=0}^2i^{3-\alpha}\int_{|y|> n^{-4}}(y^3\sigma(y))^{(\alpha)}\left(\varphi^{(3-\alpha)}(b)s_n(b+iy)
-\varphi^{(3-\alpha)}(a)s_n(a+iy)\right)dy\nonumber\\
&&+\frac{n}{2\pi}\frac{1}{3!}\int_a^b\int_{|y|> n^{-5}}\varphi'(x)(y^3\sigma(y))^{(3)}s_n(x+iy)dxdy\nonumber\\
&&+\frac{1}{2\pi}\frac{1}{3!}\sum_{l=1}^n\int_a^b\int_{|y|\leq n^{-5}}\varphi'(x)(y^3\sigma(y))^{(3)}\frac{1}{\lambda_l-z}dxdy+\mathcal{O}(n^{-2-3c_0})\nonumber\\
&&:=A_4+A_5+A_6+A_7+A_8+\mathcal{O}(n^{-2-3c_0}). \label{3.1}
\end{eqnarray}

To use the above representation in the proof of Theorem \ref{th.3.1}, we shall provide a more easily handled condition on the Stieltjes transform instead of $N_n(I)=0$. Such a trick is from Tao and Vu \cite{TV2}.
\begin{lem}\label{lem.1} For some positive constant $A_0$ (independent of $c_0$ and $a$), if
\begin{eqnarray}
Im s_n(a+in^{-1-2A_0c_0})\leq n^{-A_0c_0}, \label{5.5}
\end{eqnarray}
one has
\begin{eqnarray}
\inf_{x\in I}\min_l|\lambda_l-x|>n^{-1-A_0c_0} \label{5.8}
\end{eqnarray}
when $c_0$ is sufficiently small.
\end{lem}

A direct consequence of Lemma \ref{lem.1} is $N_n(I)=0$ when (\ref{5.5}) holds and $c_1>A_0c_0$.
\begin{proof}
By definition, one has
\begin{eqnarray*}
Im s_n(a+in^{-1-2A_0c_0})=n^{-2-2A_0c_0}\sum_{l=1}^n\frac{1}{(\lambda_l-a)^2+n^{-2-4A_0c_0}}.
\end{eqnarray*}
By the assumption $Im s_n(a+in^{-1-2A_0c_0})\leq n^{-A_0c_0}$, we can get
\begin{eqnarray*}
\min_l(\lambda_l-a)^2\geq n^{-2-A_0c_0}-n^{-2-4A_0c_0}.
\end{eqnarray*}
Consequently, when $n$ is sufficiently large one has
\begin{eqnarray}
\min_l|\lambda_l-a|\geq  2n^{-1-A_0c_0}.\label{5.6}
\end{eqnarray}
Now if $c_0$ is sufficiently small such that $c_1\geq A_0c_0$, we can easily get by triangular inequality that
\begin{eqnarray*}
\min_{x\in I}\min_l|\lambda_l-x|\geq n^{-1-A_0c_0}.
\end{eqnarray*}
Here $c_1$ is the constant in ${(\bf{b})}$ of the definition of $\mathcal{F}_n$. Thus we conclude the proof.
\end{proof}

Moreover, we have the following lemma due to Tao and Vu \cite{TV2}.
\begin{lem}\label{lem.2}\emph{(Corollary 15, \cite{TV2})} For any $a\in [-2+\delta,2-\delta]$, there exists a sufficiently large constant $A_0>0$ (independent of $c_0$ and $a$),
\begin{eqnarray*}
Im s_n(a+in^{-1-2A_0c_0})\leq n^{-A_0c_0}/2
\end{eqnarray*}
holds with high probability.
\end{lem}
\begin{rem}
The proof of the above lemma in \cite{TV2} is based on the level repulsion estimate of Wigner matrices (see Proposition 14,\cite{TV2}). The proof of the level repulsion in \cite{TV2} needs the conditions that $a\in[-2+\delta,2-\delta]$ and the distributions of the matrix elements are supported on at least three points. That is why we make these assumptions in our main results. However, we believe these restrictions are not necessary and can be removed. We will not pursue this direction in this paper.
\end{rem}

Pursuing the argument in \cite{TV2}, we define a smooth cutoff function $\chi(x)$ to the region $|x|\leq n^{-A_0c_0}$ that equals $1$ for $|x|\leq n^{-A_0c_0}/2$. Thus by Lemma \ref{lem.2}, one sees that $\chi(Im s_n(a+in^{-1-2A_0c_0}))$ is equal to $1$ with high probability. Consequently, it suffices to prove the fact that the quantity
\begin{eqnarray*}
\mathbb{E}\{G(\mathcal{L}_n[\varphi])\chi(Im s_n(a+in^{-1-2A_0c_0}))\}
\end{eqnarray*}
is highly insensitive for $\varphi\in \mathcal{F}_n$.

Moreover, we have mentioned above that by Lemma \ref{lem.1}, one has $\chi(Im s_n(a+in^{-1-2A_0c_0}))\neq 0$ implies that $N_n(I)=0$. Therefore, we can use the representation (\ref{3.1}) for $\mathcal{L}_n[\varphi]$. Moreover, by the bound on the derivative of $G$ (see (\ref{5.3})), one has
\begin{eqnarray}
&&\mathbb{E}\{G(\mathcal{L}_n[\varphi])\chi(Im s_n(a+in^{-1-2A_0c_0}))\}\nonumber\\
&&=\mathbb{E}\{G(\sum_{l=4}^8A_l)\chi(Im s_n(a+in^{-1-2A_0c_0}))\}+\mathcal{O}(n^{-2-2c_0}). \label{3.3}
\end{eqnarray}

 To show that the above quantity is highly insensitive, the main task is to provide the stability of $s_n(x+iy)$ involved in $A_i$ and $\chi(Im s_n(a+in^{-1-2A_0c_0}))$ in the swapping procedure. To this end, we need the Taylor expansion for $s_n(x+iy)$ proved by Tao and Vu \cite{TV2}: Lemma \ref{lem.3} stated in Section 2. That is to say, for Wigner matrices $M_n$ and $M'_n$, we can start from $M_n$, and then replace its elements one pair (or one unit for the diagonal case) a time by the corresponding one of  $M'_n$ and study the stability of $s_n(x+iy)$ under such a swapping process. To achieve this aim, we let $M_n^{(1)}, M_n^{(2)}$ be two adjacent matrices in the swapping procedure in the sense that we can write
 \begin{eqnarray*}
 M^{(1)}=M_0+\frac{1}{\sqrt{n}}\xi^{(1)}V,\quad M^{(2)}=M_0+\frac{1}{\sqrt{n}}\xi^{(2)}V
 \end{eqnarray*}
 for some elementary matrix $V$. And $\xi^{(1)}, \xi^{(2)}$ are two real random variables matching to the fourth order and bounded in magnitude by $n^{\mathcal{O}(c_0)}$. Moreover, $M_0$ is independent of $\xi^{(1)}$ and $\xi^{(2)}$.

 To describe the swapping process, we use the notation $s_{\xi^{(1)}}(x+iy)$ to denote the Stieltjes transform for $M_n^{(1)}$, and $s_{\xi^{(2)}}(x+iy)$ for $M_n^{(2)}$. Correspondingly we distinguish $\mathcal{L}_n[\varphi]$ by $\mathcal{L}_n^{(1)}[\varphi]$ and $\mathcal{L}_n^{(2)}[\varphi]$ for $M_n^{(1)}$ and $M_n^{(2)}$ respectively. Besides, we use the notation $A_l^{(1)}$ and $A_l^{(2)}$ to denote $A_l$ ($l=4,\cdots,8$ ) for $M_n^{(1)}$ and $M_n^{(2)}$ respectively. With these notations, we will show the quantity
 \begin{eqnarray*}
\mathbb{E}\{G(\sum_{l=4}^8A_l^{(1)})\chi(Im s_{\xi^{(1)}}(a+in^{-1-2A_0c_0}))\}
 \end{eqnarray*}
 only changes by $\mathcal{O}(n^{-2-\mathcal{O}(c_0)})$ when $\xi^{(1)}$ is replaced by $\xi^{(2)}$ in the off diagonal case, or $\mathcal{O}(n^{-1-\mathcal{O}(c_0)})$ in the diagonal case. Then by a telescoping arguments, after $\mathcal{O}(n^2)$ steps of replacement, we can easily get that
\begin{eqnarray*}
\mathbb{E}\{G(\mathcal{L}_n[\varphi])\chi(Im s_n(a+in^{-1-2A_0c_0}))\}
\end{eqnarray*}
is highly insensitive.

 In order to apply Lemma \ref{lem.3}, one shall guarantee the condition (\ref{3.2}). We need the following crucial lemma.
 \begin{lem} \label{lem.4}(Uniform resolvent bounds). We have the following two assertions on the resolvent bounds.

 (i): If $\chi(Im s_{\xi^{(1)}}(a+in^{-1-2A_0c_0}))\neq 0$, then with overwhelming probability
 \begin{eqnarray}
 \sup_{x\in I}\sup_{y>0}||R_{\xi^{(1)}}(x+iy)||_{(\infty,1)}= \mathcal{O}(n^{\mathcal{O}(c_0)}) \label{3.7}
 \end{eqnarray}
 and
 \begin{eqnarray}
 \sup_{x\in I}\sup_{y>0}||R_{0}(x+iy)||_{(\infty,1)}= \mathcal{O}(n^{\mathcal{O}(c_0)}).\label{3.8}
 \end{eqnarray}

 (ii): If $y_0=n^{-1-c_0}$, then with overwhelming probability
 \begin{eqnarray}
 \sup_{x\in\mathcal{U}_\epsilon}\sup_{|y|>y_0}||R_{\xi^{(1)}}(x+iy)||_{(\infty,1)}= \mathcal{O}(n^{\mathcal{O}(c_0)})\label{3.9}
 \end{eqnarray}
 and
 \begin{eqnarray}
 \sup_{x\in\mathcal{U}_\epsilon}\sup_{|y|>y_0}||R_{0}(x+iy)||_{(\infty,1)}= \mathcal{O}(n^{\mathcal{O}(c_0)}).\label{3.10}
 \end{eqnarray}
 \end{lem}
 \begin{proof} At first, we prove (i). We learn from the proof of Lemma \ref{lem.1} that when
 $$\chi(Im s_{\xi^{(1)}}(a+in^{-1-2A_0c_0}))\neq 0,$$
 there exists (\ref{5.6}).
 Besides, by the spectral decomposition we can easily get
 \begin{eqnarray}
 ||R_{\xi^{(1)}}(x+iy)||_{(\infty,1)}\leq \sum_{l=1}^{n}\frac{\sup_{1\leq j\leq n}||u_j(M_n^{(1)})||_{l^\infty}^2}{|\lambda_l(M_n^{(1)})-x-iy|},\label{3.5}
 \end{eqnarray}
 where $u_j(M_n^{(1)})$ is the unit eigenvector of $M_n^{(1)}$ corresponding to $\lambda_j(M_n^{(1)})$. With (\ref{5.1}) in Lemma \ref{lem.5.1}, one has
 \begin{eqnarray*}
 ||R_{\xi^{(1)}}(x+iy)||_{(\infty,1)}\leq n^{-1+\mathcal{O}(c_0)}\sum_{l=1}^n\frac{1}{|\lambda_l(M_n^{(1)})-x|}.
 \end{eqnarray*}
 For $x\in I$, by using (\ref{5.6}) we have
 \begin{eqnarray*}
 |\lambda_l(M_n^{(1)})-x|\geq|\lambda_l(M_n^{(1)})-a|-n^{-1-c_1}\geq  \frac12|\lambda_l(M_n^{(1)})-a|
 \end{eqnarray*}
 when $c_1>A_0c_0$. Consequently we have
 \begin{eqnarray*}
 \sup_{x\in I}\sup_{y>0}||R_{\xi^{(1)}}(x+iy)||_{(\infty,1)}\leq 2n^{-1+\mathcal{O}(c_0)}\sum_{l=1}^n\frac{1}{|\lambda_l(M_n^{(1)})-a|}.
 \end{eqnarray*}
 By the argument in Tao and Vu (the proof of Lemma 16, \cite{TV2}), if $\chi(Im s_{\xi^{(1)}}(a+in^{-1-2A_0c_0}))\neq 0$, one has with overwhelming probability that
 \begin{eqnarray*}
 \sum_{l=1}^n\frac{1}{|\lambda_l(M_n^{(1)})-a|}=\mathcal{O}(n^{1+\mathcal{O}(c_0)}).
 \end{eqnarray*}
 Thus we have
 \begin{eqnarray*}
 \sup_{x\in I}\sup_{y>0}||R_{\xi^{(1)}}(x+iy)||_{(\infty,1)}= \mathcal{O}(n^{\mathcal{O}(c_0)}).
 \end{eqnarray*}
 For $R_0(x+iy)$, we use the fact that for $y>0$, when $|t|||R_t||_{(\infty,1)}=o(\sqrt{n})$,
 \begin{eqnarray}
 ||R_0||_{(\infty,1)}\leq (1+o(1))||R_t||_{(\infty,1)}. \label{3.6}
 \end{eqnarray}
 (\ref{3.6}) is a consequence of Neumann series formula, we refer to \cite{TV2} for the details of the proof.
 Thus we also have
 \begin{eqnarray*}
 \sup_{x\in I}\sup_{y>0}||R_{0}(x+iy)||_{(\infty,1)}= \mathcal{O}(n^{\mathcal{O}(c_0)}).
 \end{eqnarray*}
 Now we turn to the proof of (ii). For a sufficiently large constant $A$, we set $\eta=n^{-1+Ac_0}$. We cover the interval $\mathcal{U}_\epsilon$ by the union of the intervals $I_k=[(k-\frac12)\eta,(k+\frac12)\eta]$ with the integer index $k$ running from $-\lfloor(2+\epsilon)\eta^{-1}\rfloor-1$ to $\lfloor(2+\epsilon)\eta^{-1}\rfloor+1$.

 Now note that by (\ref{3.5}) we have
 \begin{eqnarray*}
 ||R_{\xi^{(1)}}(x+iy)||_{(\infty,1)}&\leq& n^{-1+\mathcal{O}(c_0)}\sum_{l=1}^{n}\frac{1}{\max\{|\lambda_l(M_n^{(1)})-x|,|y|\}}\nonumber\\
 &=&n^{-1+\mathcal{O}(c_0)}\sum_{k}\sum_{l:\lambda_l(M_n^{(1)})\in I_k}\frac{1}{\max\{|\lambda_l(M_n^{(1)})-x|,|y|\}}.
 \end{eqnarray*}
 By the fact that $N_J=O(n|J|)$ with overwhelming probability for any interval $J$ with length $|J|\geq n^{-1+Ac_0}$ (see Lemma \ref{lem.5.1}), we can immediately get
 \begin{eqnarray*}
 \sup_{x\in\mathcal{U}_\epsilon}\sup_{y>y_0}||R_{\xi^{(1)}}(x+iy)||_{(\infty,1)}=\mathcal{O}(n^{\mathcal{O}(c_0)})
 \end{eqnarray*}
with overwhelming probability. Again by the resolvent bound (\ref{3.6}) one has
 \begin{eqnarray*}
 \sup_{x\in\mathcal{U}_\epsilon}\sup_{y>y_0}||R_{0}(x+iy)||_{(\infty,1)}=\mathcal{O}(n^{\mathcal{O}(c_0)}).
 \end{eqnarray*}
 Thus we complete the proof.
 \end{proof}

 If we condition on the event that (\ref{3.8}) and (\ref{3.10}) hold, then we also have (\ref{3.7}) and (\ref{3.9}) by swapping the roles of $R_0$ and $R_t$ in (\ref{3.6}). Since $|\xi^{(1)}|\leq n^{\mathcal{O}(c_0)}$, together with (\ref{3.8}) and (\ref{3.10}), we can use Lemma \ref{lem.3}. Besides, under the condition $\chi(Im s_{\xi^{(1)}}(a+in^{-1-2A_0c_0}))\neq 0$, we have (\ref{5.8}), which trivially implies
 \begin{eqnarray*}
 \sup_{x\in I}\sup_y |\frac1n\sum_{l=1}^n\frac{1}{\lambda_l-z}|=\mathcal{O}(n^{1+\mathcal{O}(c_0)}).
 \end{eqnarray*}
 Consequently, we have
 \begin{eqnarray*}
 A_8^{(1)}=\mathcal{O}(n^{-3+\mathcal{O}(c_0)}).
\end{eqnarray*}
Now we set
\begin{eqnarray*}
&&A_4^0=\frac{n}{2\pi}\int_{-2-\epsilon}^{2+\epsilon}\left(\int_{-1}^{-1/2}
+\int_{1/2}^{1}\right)\left(\sum_{j=0}^{3}\frac{(-1)^j}{j!}i(y^j\sigma'(y))^{(j)}\right)\varphi(x)s_0(x+iy)dxdy,\nonumber\\
&&A_5^0=\frac{n}{2\pi}\left(\int_{b}^{2+\epsilon}+\int_{-2-\epsilon}^{a}\right)\int_{|y|> y_0}\frac{1}{3!}\varphi^{(4)}(x)(iy)^3\sigma(y)s_0(x+iy)dxdy,\nonumber\\
&&A_6^0=\frac{n}{2\pi}\frac{1}{3!}\sum_{\alpha=0}^2i^{3-\alpha}\int_{|y|>n^{-4}}(y^3\sigma(y))^{(\alpha)}\left(\varphi^{(3-\alpha)}(b)s_0(b+iy)
-\varphi^{(3-\alpha)}(a)s_0(a+iy)\right)dy,\nonumber\\
&&A_7^0=\frac{n}{2\pi}\frac{1}{3!}\int_a^b\int_{|y|> n^{-5}}\varphi'(x)(y^3\sigma(y))^{(3)}s_0(x+iy)dxdy.
\end{eqnarray*}
Then we let
\begin{eqnarray*}
&&G^{(1)}_n=:G(\sum_{l=4}^7A_l^{(1)})\chi(Im s_{\xi^{(1)}}(a+in^{-1-2A_0c_0})),\nonumber\\
&&G^{0}_n=:G(\sum_{l=4}^7A_l^0)\chi(Im s_0(a+in^{-1-2A_0c_0})).
 \end{eqnarray*}
 Correspondingly we can define $G_n^{(2)}$. Our aim is to expand $G^{(1)}_n$ around $G^{(0)}_n$. We formulate the result as the following lemma.
 \begin{lem} \label{lem.5}With the above notations, when we condition on the event that (\ref{3.8}) and (\ref{3.10}) hold, we have
 \begin{eqnarray}
 G^{(1)}_n=G^{0}_n+\sum_{j=1}^4d_j(\xi^{(1)})^j+\mathcal{O}(n^{-5/2+\mathcal{O}(c_0)}), \label{3.11}
 \end{eqnarray}
  where the coefficients $d_j$ are independent of $\xi^{(1)}$ and obey the bounds
 \begin{eqnarray*}
 d_j=\mathcal{O}(n^{-j/2+\mathcal{O}(c_0)}), \quad j=1,\cdots.4.
 \end{eqnarray*}
 \end{lem}
 \begin{proof}
 Note that if (\ref{3.8}) and (\ref{3.10}) hold, we can use Lemma \ref{lem.3} to expand $s_{\xi^{(1)}}$ (resp. $A_l^{(1)}$) around $s_0$ (resp. $A_l^{0}$).  Then by the assumption on the derivatives of $G(x)$ and the fact that
 \begin{eqnarray*}
 \frac{d^j}{dx^j}\chi(x)=\mathcal{O}(n^{\mathcal{O}(jc_0)}),
 \end{eqnarray*}
 we can conclude (\ref{3.11}) by applying Taylor expansion to $G(\cdot)$ and $\chi(\cdot)$.
 \end{proof}

 \begin{proof}[Proof of Theorem \ref{th.3.1}] Using Lemma \ref{lem.5} to both $M^{(1)}$ and $M^{(2)}$, we have
  \begin{eqnarray*}
  G^{(1)}_n=G^{0}_n+\sum_{j=1}^4d_j(\xi^{(1)})^j+\mathcal{O}(n^{-5/2+\mathcal{O}(c_0)})
  \end{eqnarray*}
  and
  \begin{eqnarray*}
  G^{(2)}_n=G^{0}_n+\sum_{j=1}^4d_j(\xi^{(2)})^j+\mathcal{O}(n^{-5/2+\mathcal{O}(c_0)}).
  \end{eqnarray*}
  Then taking expectation with respect to $\xi^{(1)}$ and $\xi^{(2)}$ respectively, and by the telescope arguments on $\mathcal{O}(n^2)$ steps of swapping, we can immediately get the conclusion by the matching moments assumption.
 \end{proof}

 As an application of Theorem \ref{th.3.1}, we can prove our main result Theorem \ref{th.2}.
 \begin{proof}[Proof of Theorem \ref{th.2}] Note that Lemma \ref{lem.11} holds for Wigner matrices under the conditions of Theorem \ref{th.2}.
By the argument in Section 3, it suffices to show that (ii) and (iii) of Theorem \ref{th.1} are still valid for general Wigner matrices. To combine Theorem \ref{th.1} with Theorem \ref{th.3.1}, we let $M_n$ be a general Wigner matrix and $M'_n$ be GUE. Moreover, we shall truncate the elements of $W_n$ and $W'_n$ at $\mathcal{O}(n^{c_0})$ to adapt to the condition of Theorem \ref{th.3.1}. Under the condition $\mathbf{C_0}$, it is easy to see such a truncation does not alter the limiting behavior of both two types of partial linear eigenvalue statistics. Moreover, the truncation will change the first four moments of the elements by only $\mathcal{O}(e^{-n^{\mathcal{O}(c_0)}})$, which can be absorbed in the remainder $\mathcal{O}(n^{-5/2+\mathcal{O}(c_0)})$ when we take expectations on both sides of (\ref{3.11}). Next, we define a smooth modification of $f_u$ by
\begin{eqnarray}
h_u(x)=:h_u(n,x)=(f(x)-f(u))\tilde{\chi}_u(n,x), \label{5.2}
\end{eqnarray}
where $\tilde{\chi}_u(n,x)$ is an $n$-dependent smooth cutoff to the region $x\in (-\infty, u+n^{-1-c_1}/2)$ that equals 1 for $x\in (-\infty, u-n^{-1-c_1}/2)$. Similar to (\ref{2.8}), one has
\begin{eqnarray*}
\mathcal{A}_n[f_u;u]=\mathcal{L}_n[f_u]=\mathcal{L}_n[h_u]+o(1)
\end{eqnarray*}
with overwhelming probability for all Wigner matrices satisfying {\emph{Condition}} $\mathbf{C_0}$.

  It is not difficult to see that to extend (ii) in Theorem \ref{th.1}, it suffices to show that for any interval $J=[j_1,j_2]$,
 \begin{eqnarray}
 \mathbb{P}(\mathcal{L}_n^{M'}[h_u]\in J_{-})-n^{-c_0}\leq \mathbb{P}(\mathcal{L}_n^M[h_u]\in J)\leq \mathbb{P}(\mathcal{L}_n^{M'}[h_u]\in J_{+})+n^{-c_0}. \label{3.12}
 \end{eqnarray}
 Here
 \begin{eqnarray*}
 J_+=[j_1-n^{-c_0/10}, j_2+n^{-c_0/10}],\quad J_-=[j_1+n^{-c_0/10}, j_2-n^{-c_0/10}].
 \end{eqnarray*}
 We only show the second inequality of (\ref{3.12}) below since the first one is analogous. Let $G:\mathbb{R}\rightarrow \mathbb{R}^+$ equal to one in $J$ and vanish outside of $J_+$ such that
 \begin{eqnarray*}
 \frac{d^j}{dx^j}G(x)=\mathcal{O}(n^{c_0}),\quad j=0,\cdots,5.
 \end{eqnarray*}
 We can apply Theorem \ref{th.3.1} to $G$ defined above. Observe that
 \begin{eqnarray*}
 \mathbb{P}(\mathcal{L}_n^M[h_u]\in J)\leq \mathbb{E}G(\mathcal{L}_n^{M}[h_u]),
 \end{eqnarray*}
 and
 \begin{eqnarray*}
 \mathbb{E}G(\mathcal{L}_n^{M'}[h_u]) \leq  \mathbb{P}(\mathcal{L}_n^{M'}[h_u]\in J_+).
 \end{eqnarray*}
 Thus by Theorem \ref{th.3.1}, we immediately get the second inequality of (\ref{3.12}).  Thus we have proved that (ii) of Theorem \ref{th.1} is still valid for general Wigner matrices. For (iii), by using (\ref{2.18}) and (\ref{2.19}) again, we can also reduce the problem to compare $\mathcal{L}_n^M[f_{\gamma_{k/n}}]$ and $\mathcal{L}_n^{M'}[f_{\gamma_{k/n}}]$. Then we only need to replace $u$ by $\gamma_{k/n}$ in the arguments above and get the conclusion. Thus we finally complete the proof of Theorem \ref{th.2}.
 \end{proof}
 \section{Partial sum process for Wigner matrices}
 In this section, we prove Theorem \ref{th.4} by providing the finite-dimensional convergence and tightness of the sequence $\mathcal{S}_n[f;t]$. Our strategy is to use the comparison theorem established in the last section to both two parts of the proof. Below we will use the notation $\mathcal{S}^M_n[f;t]$ and $\mathcal{S}^{M'}_n[f;t]$ to denote $\mathcal{S}_n[f;t]$ for Wigner matrices $M_n$ and $M'_n$ respectively. Moreover, we will specify $M'_n$ to be GUE in this section. At first, we will prove the following lemma on the expectations.
 \begin{lem} \label{lem.6.1}Under the assumption of Theorem \ref{th.4}, for $t\in[\delta,1-\delta]$, we have
 \begin{eqnarray*}
 \mathbb{E}\mathcal{S}^M_n[f;t]-\mathbb{E}\mathcal{S}^{M'}_n[f;t]=\mathcal{O}(n^{-c})
 \end{eqnarray*}
 for some fixed constant $c>0$.
 \end{lem}
 \begin{proof}
By the assumptions on $f(x)$ in Theorem \ref{th.4} and the large deviation estimate of extreme eigenvalue in Lemma \ref{lem.12}, we can and do assume that $f(x)$ is compactly supported on $\mathcal{U}_\epsilon$ below. By definition, we shall provide that
 \begin{eqnarray}
 \mathbb{E}\mathcal{B}^M_n[f;\lfloor nt\rfloor]-\mathbb{E}\mathcal{B}^{M'}_n[f;\lfloor nt\rfloor]=\mathcal{O}(n^{-c}), \label{6.1}
 \end{eqnarray}
 and
 \begin{eqnarray}
 \mathbb{E}f(\lambda_{\lfloor nt\rfloor+1}(M_n))-\mathbb{E}f(\lambda_{\lfloor nt\rfloor+1}(M'_n))=\mathcal{O}(n^{-c}). \label{6.2}
 \end{eqnarray}
 Relying on (\ref{2.18}) and the fact that both $\mathcal{B}^M_n[f;\lfloor nt\rfloor]$ and $\mathcal{B}^{M'}_n[f;\lfloor nt\rfloor]$ are $\mathcal{O}(n)$, one can prove
 \begin{eqnarray*}
  \mathbb{E}\widehat{\mathcal{B}}^M_n[f;\lfloor nt\rfloor]-\mathbb{E}\widehat{\mathcal{B}}^{M'}_n[f;\lfloor nt\rfloor]=\mathcal{O}(n^{-c})
 \end{eqnarray*}
 instead of (\ref{6.1}). Moreover, by (\ref{2.19}) and the definition of $h_u(x)$ in (\ref{5.2}), apparently it suffices to show
 \begin{eqnarray}
 \mathbb{E}\mathcal{L}^M_n[h_{u(t)}]-\mathbb{E}\mathcal{L}^{M'}_n[h_{u(t)}]=\mathcal{O}(n^{-c}). \label{6.3}
 \end{eqnarray}
 Here we recall the notation $u(t)=\gamma_{\lfloor nt\rfloor/n}$.

Let $\chi_3(x)$ be a smooth cutoff function to the region $[-2n^{c_0/2}, 2n^{c_0/2}]$ which equals to $1$ in $[-n^{c_0/2},n^{c_0/2}]$. Now we define the function $G_1$ as
\begin{eqnarray*}
G_1(x)=(x-\sum_{l=1}^n h_{u(t)}(\gamma_{l/n}))\chi_3(x-\sum_{l=1}^n h_{u(t)}(\gamma_{l/n})).
\end{eqnarray*}
Obviously, $G_1(x)$ satisfies the condition (\ref{5.3}).

Now we claim that
\begin{eqnarray}
\mathbb{E}G_1(\mathcal{L}_n[h_{u(t)}])=\mathbb{E}\mathcal{L}_n[h_{u(t)}]-\sum_{l=1}^n h_{u(t)}(\gamma_{l/n})+\mathcal{O}(n^{-c}). \label{6.4}
\end{eqnarray}
To see (\ref{6.4}), it suffices to show that
\begin{eqnarray}
\chi_3(\mathcal{L}_n[h_{u(t)}]-\sum_{l=1}^n h_{u(t)}(\gamma_{l/n}))=1 \label{6.5}
\end{eqnarray}
holding with overwhelming probability. To this end, we use the rigidity property stated in Lemma \ref{lem.9}. With the aid of this lemma, we now show the validity of (\ref{6.5}) as follows.
Note that with overwhelming probability one has
\begin{eqnarray*}
&&|\mathcal{L}_n[h_{u(t)}]-\sum_{l=1}^n h_{u(t)}(\gamma_{l/n})|
\leq \sup_{x\in\mathcal{U}_\epsilon}|h'_{u(t)}(x)|\sum_{l=1}^n|\lambda_l-\gamma_{l/n}|\nonumber\\
&&\leq C(\log n)^{C\log\log n}\sum_{l=1}^n [\min(l,n-l+1)]^{-1/3}n^{-2/3}
\leq n^{c_0/2}
\end{eqnarray*}
for $n$ sufficiently large. Thus we have shown (\ref{6.4}).

Consequently, it suffices to show (\ref{6.2}) and
\begin{eqnarray}
\mathbb{E}G_1(\mathcal{L}^M_n[h_{u(t)}])-\mathbb{E}G_1(\mathcal{L}^{M'}_n[h_{u(t)}])=\mathcal{O}(n^{-c}). \label{6.6}
\end{eqnarray}
Observe that (\ref{6.6}) is a direct consequence of Theorem \ref{th.3.1} (After a harmless truncation towards the elements of both $W_n$ and $W'_n$). Thus we only need to show (\ref{6.2}) in the sequel.
Note that
\begin{eqnarray*}
&&|f(\lambda_{\lfloor nt\rfloor+1}(M_n))-f(\lambda_{\lfloor nt\rfloor+1}(M'_n))|\nonumber\\
&&\leq \sup_{x\in \mathcal{U}_\epsilon}|f'(x)||\lambda_{\lfloor nt\rfloor+1}(M_n)-\lambda_{\lfloor nt\rfloor+1}(M'_n)|\nonumber\\
&&\leq C\frac{(\log n)^{C\log\log n}}{n}
\end{eqnarray*}
holding with overwhelming probability. By the assumption that $f(x)$ is compactly supported thus bounded, one can immediately get (\ref{6.2}). Hence, we complete the proof.
 \end{proof}
 Now we begin to prove the finite dimensional convergence of the sequence $(\mathcal{S}_n^M[f;t])^\circ$. We formulate the result as the following lemma and then prove it.
 \begin{lem} \label{lem.6.3}Under the assumption of Theorem \ref{th.4}, for any fixed positive integer $r$ and points $t_1,\cdots,t_r\in[\delta,1-\delta]$, and for any fixed numbers $\alpha_1,\cdots,\alpha_r\in \mathbb{R}$, we have
\begin{eqnarray*}
\sum_{l=1}^r\alpha_l(\mathcal{S}^M_n[f;t_l])^\circ\stackrel{d}\longrightarrow \sum_{l=1}^r\alpha_l\mathcal{S}[f;t_l].
\end{eqnarray*}
 \end{lem}
 \begin{proof}
 By the same discussion as that for GUE taken in Section 4, we can transfer the problem to show
 \begin{eqnarray*}
 (\mathcal{L}^M_n[\sum_{l=1}^r\alpha_lf_{\gamma_{t_l}}])^\circ \stackrel{d}\longrightarrow \sum_{l=1}^r\alpha_l\mathcal{S}[f;t_l]
 \end{eqnarray*}
 instead. As mentioned in Section 4, $\sum_{l=1}^r\alpha_lf_{\gamma_{t_l}}$ is a continuous function with $r$ possibly non-differentiable points $t_1,\cdots,t_r$. Now we choose $r$ intervals $J_1,\cdots, J_r$ containing $\gamma_{t_1},\cdots,\gamma_{t_r}$ respectively with lengths $|J_l|\leq n^{-1-c_1}$ for all $l=1,\cdots,r$. And then we define a smooth modification function $h_{t_1,\cdots, t_r}$ which coincides with $\sum_{l=1}^r\alpha_lf_{\gamma_{t_l}}$ on $\mathbb{R}\setminus\cup_{l=1}^rJ_l$. Thus, it is easy to see $h_{t_1,\cdots, t_r}\in \mathcal{F}_n^r$. Now by using Theorem \ref{th.3.1} and a routine discussion as that in the proof of Theorem \ref{th.2}, one can get that
 \begin{eqnarray}
 \mathcal{L}^M_n[h_{t_1,\cdots,t_r}]-\mathbb{E}\mathcal{L}_n^{M'}[h_{t_1,\cdots,t_r}]\stackrel{d}\longrightarrow \sum_{l=1}^r\alpha_l\mathcal{S}[f;t_l]. \label{6.7}
 \end{eqnarray}
 Here $M'=M'_n$ is GUE. Apparently, by (\ref{6.7}) one can get
 \begin{eqnarray*}
 \mathcal{L}^M_n[\sum_{l=1}^r\alpha_lf_{\gamma_{t_l}}]-\mathbb{E}\mathcal{L}^{M'}_n[\sum_{l=1}^r\alpha_lf_{\gamma_{t_l}}] \stackrel{d}\longrightarrow \sum_{l=1}^r\alpha_l\mathcal{S}[f;t_l].
 \end{eqnarray*}
 Then by Lemma \ref{lem.6.1}, we can conclude the proof of Lemma \ref{lem.6.3}.
 \end{proof}

 Thus the only thing left is to prove the tightness of the sequence $(\mathcal{S}_n^M[f;t])^\circ$. It suffices to provide the H\"{o}lder condition: for arbitrary $t,s \in [\delta,1-\delta]$ and sufficiently large $n$
 \begin{eqnarray*}
 \mathbb{E}|(\mathcal{S}_n^M[f;t])^\circ-(\mathcal{S}_n^M[f;s])^\circ|^2\leq C|t-s|^\alpha
 \end{eqnarray*}
 holds for some constants $C>0$, and $\alpha>1$ both independent of $t,s$.

 At first, we claim that under the assumptions in Theorem \ref{th.4}, we always have for $l\in [\delta n, (1-\delta)n]$
 \begin{eqnarray}
Var\{f(\lambda_l)\}\leq C\frac{\log n}{n^2} \label{6.12}
 \end{eqnarray}
with some positive constant $C$ depending only on $\delta$ and $f$. Based on Lemma \ref{lem.10}, the proof of (\ref{6.12}) is analogous to that of (\ref{4.1}). The only difference is that we need to use Lemma \ref{lem.12} instead of (\ref{2.33}). Here we omit the detail.

 Similar to the GUE case, we always assume that $s\leq t$ and separate the issue into three cases: $t-s\leq n^{-1}$, $n^{-1}\leq t-s\leq \eta_n$ and $t-s\geq \eta_n$. It is easy to see the discussions for the first two cases for GUE in Section 4 are also valid for general Wigner matrix with the aid of (\ref{6.12}). Thus we will focus on the third case in the sequel.
 \begin{lem} \label{lem.6.4} Under the assumptions of Theorem \ref{th.4}, for $t,s\in[\delta,1-\delta]$ such that $t-s\geq \eta_n=:\log^{-L}n$ with sufficiently large constant $L$, one has for $n$ large enough
 \begin{eqnarray*}
 \mathbb{E}|(\mathcal{S}_n^M[f;t])^\circ-(\mathcal{S}_n^M[f;s])^\circ|^2\leq C|t-s|^\alpha
 \end{eqnarray*}
for some constants $C>0$, and $\alpha>1$ both independent of $t,s$.
 \end{lem}
 \begin{proof}
 To prove Lemma \ref{lem.6.4}, we will use Theorem \ref{th.3.1} to compare the general case with the Gaussian case. At first, note that it suffices to prove
 \begin{eqnarray*}
  \mathbb{E}|(\mathcal{B}_n^M[f;\lfloor nt\rfloor])^\circ-(\mathcal{B}_n^M[f;\lfloor ns\rfloor])^\circ|^2\leq C|t-s|^\alpha
 \end{eqnarray*}
 instead since
 \begin{eqnarray*}
Var f(\lambda_{\lfloor nt\rfloor+1})+Var f(\lambda_{\lfloor ns\rfloor+1})\leq C(\delta)\frac{\log^2n}{n}.
 \end{eqnarray*}
 Moreover, by (2.18) we have
 \begin{eqnarray}
 &&\mathbb{E}|(\mathcal{B}_n^M[f;\lfloor nt\rfloor])^\circ-(\mathcal{B}_n^M[f;\lfloor ns\rfloor])^\circ|^2\nonumber\\
 &&=\mathbb{E}|(\widehat{\mathcal{B}}_n^M[f;\lfloor nt\rfloor])^\circ-(\widehat{\mathcal{B}}_n^M[f;\lfloor ns\rfloor])^\circ|^2+\mathcal{O}(\frac{(\log n)^{C\log\log n}}{n})\nonumber\\
 &&=Var\{\mathcal{L}^M_n[f_{u(t)}-f_{u(s)}]\}+\mathcal{O}(\frac{(\log n)^{C\log\log n}}{n})\nonumber\\
 &&=Var\{\mathcal{L}^M_n[h_{u(t)}-h_{u(s)}]\}+\mathcal{O}(n^{-c}). \label{6.8}
 \end{eqnarray}
 Here $h_{u}$ is defined in (\ref{5.2}).

 Note that we have (\ref{2.24}) for GUE. Thus by (\ref{6.8}), we see that it suffices to show
 \begin{eqnarray}
 Var\{\mathcal{L}^M_n[h_{u(t)}-h_{u(s)}]\}-Var\{\mathcal{L}^{M'}_n[h_{u(t)}-h_{u(s)}]\}=\mathcal{O}(n^{-c}). \label{6.9}
 \end{eqnarray}
Observe that
 \begin{eqnarray*}
 Var\{\mathcal{L}^M_n[h_{u(t)}-h_{u(s)}]\}&=&\mathbb{E}(\mathcal{L}_n^M[h_{u(t)}-h_{u(s)}]-\mathbb{E}\mathcal{L}_n^{M'}[h_{u(t)}-h_{u(s)}])^2\nonumber\\
 &&-(\mathbb{E}\{\mathcal{L}_n^M[h_{u(t)}-h_{u(s)}]\}-\mathbb{E}\{\mathcal{L}_n^{M'}[h_{u(t)}-h_{u(s)}]\})^2\nonumber\\
 &=&\mathbb{E}(\mathcal{L}_n^M[h_{u(t)}-h_{u(s)}]-\mathbb{E}\mathcal{L}_n^{M'}[h_{u(t)}-h_{u(s)}])^2+\mathcal{O}(n^{-2c}).
 \end{eqnarray*}
 Here the last step follows from (\ref{6.6}).
 Thus it suffices to show that
 \begin{eqnarray*}
 \mathbb{E}(\mathcal{L}_n^M[h_{u(t)}-h_{u(s)}]-\mathbb{E}\mathcal{L}_n^{M'}[h_{u(t)}-h_{u(s)}])^2
 -Var\{\mathcal{L}_n^{M'}[h_{u(t)}-h_{u(s)}]\}=\mathcal{O}(n^{-2c}).
 \end{eqnarray*}
 Now we set
 \begin{eqnarray*}
 G_2(x)=(x-\mathbb{E}\mathcal{L}_n^{M'}[h_{u(t)}-h_{u(s)}])\chi_3(x-\mathbb{E}\mathcal{L}_n^{M'}[h_{u(t)}-h_{u(s)}]).
 \end{eqnarray*}

 Note that again by the rigidity property stated in Lemma \ref{lem.9}, one has
 \begin{eqnarray*}
 |\mathcal{L}_n^M[h_{u(t)}-h_{u(s)}]-\sum_{l=1}^n(h_{u(t)}-h_{u(s)})(\gamma_{l/n})|\leq (\log n)^{C\log\log n}
 \end{eqnarray*}
 with overwhelming probability. Simultaneously, one  can get that
 \begin{eqnarray*}
 |\mathbb{E}\mathcal{L}_n^{M'}[h_{u(t)}-h_{u(s)}]-\sum_{l=1}^n(h_{u(t)}-h_{u(s)})(\gamma_{l/n})|\leq (\log n)^{C\log\log n}.
 \end{eqnarray*}
 Thus
 \begin{eqnarray*}
 \chi_3(\mathcal{L}_n^M[h_{u(t)}-h_{u(s)}]-\mathbb{E}\mathcal{L}_n^{M'}[h_{u(t)}-h_{u(s)}])=1
 \end{eqnarray*}
 with overwhelming probability. Consequently, it suffices to compare
 \begin{eqnarray*}
 \mathbb{E}\{G_2(\mathcal{L}_n^M[h_{u(t)}-h_{u(s)}])\}
 \end{eqnarray*}
 and
  \begin{eqnarray*}
 \mathbb{E}\{G_2(\mathcal{L}_n^{M'}[h_{u(t)}-h_{u(s)}])\}.
 \end{eqnarray*}
 By using Theorem \ref{th.3.1} again, we can get
 \begin{eqnarray*}
 \mathbb{E}\{G_2(\mathcal{L}_n^M[h_{u(t)}-h_{u(s)}])\}-\mathbb{E}\{G_2(\mathcal{L}_n^{M'}[h_{u(t)}-h_{u(s)}])\}=\mathcal{O}(n^{-c}).
 \end{eqnarray*}
 Thus (\ref{6.9}) follows. Moreover, the remainder $\mathcal{O}(n^{-c})$ in (\ref{6.9}) is uniform in $t,s$, which can be seen by a careful check throughout the whole proof process. We leave it to the reader.
 \end{proof}
 Combing the finite dimensional convergence and the tightness we finally complete the proof of Theorem \ref{th.4}.
 \section{Appendix}
 In this appendix, we present some existing results on the local behavior of the spectrum of Wigner matrices satisfying {\emph{Condition}} $\mathbf{C_0}$, which can be found in the recent work on the universality property of RMT. One can refer to the series \cite{ESY}-\cite{EYY} and \cite{TV1}-\cite{TV3} for instance. We also remark here the results stated below may be proved in their original articles under weaker conditions than those made in our paper. For ease of presentation, we reformulate them under the {\emph{Condition}} $\mathbf{C_0}$ without further explanation.
 \begin{lem} \label{lem.5.1}Let $M_n$ be a Wigner matrix satisfying the Condition $\mathbf{C_0}$, one has for any interval $I$ with its length $|I|\geq n^{-1+Ac_0}$ for some sufficiently large constant $A>0$
 \begin{eqnarray}
|N_n(I)-n\int_I\rho_{sc}(x)dx|\leq \varepsilon n|I| \label{6.10}
 \end{eqnarray}
 with overwhelming probability. Moreover, one has
 \begin{eqnarray}
\sup_{1\leq j\leq n}||u_j(W_n)||_{l^\infty}\leq n^{-1/2+\mathcal{O}(c_0)} \label{5.1}
 \end{eqnarray}
 with overwhelming probability. Here $u_j(W_n)$ is the unit eigenvector corresponding to $\lambda_j(W_n)$.
 \end{lem}
 \begin{proof}
 See Theorem 1.8 and Proposition 1.10 of \cite{TV4} for instance.
 \end{proof}
 The second main lemma is an explicit description on the location of the eigenvalues proved in \cite{EYY}, named as the rigidity property for eigenvalues.
 \begin{lem}[Rigidity of eigenvalues]\label{lem.9}
Suppose that $M_n$ is a Wigner matrix obeying the Condition $\mathbf{C_0}$. One has for some positive constants $C,C',c$
\begin{eqnarray}
&&P(\exists j: \lambda_j-\gamma_{j/n}\geq (\log n)^{C\log\log n}[\min (j, n-j+1)]^{-1/3}n^{-2/3})\nonumber\\
&&\leq C'\exp[-(\log n)^{c\log\log n}] \label{6.11}
\end{eqnarray}
for $n$ large enough. Moreover, one has
\begin{eqnarray*}
\mathbb{P}(\sup_{|x|\leq 5}n|F_n(x)-F_{sc}(x)|\geq (\log n)^{C\log\log n})\leq C'\exp[-(\log n)^{c\log\log n}]
\end{eqnarray*}
for sufficiently large $n$.
\end{lem}
\begin{rem}
It is not difficult to deduce (\ref{6.10}) in Lemma \ref{lem.5.1} from (\ref{6.11}) indeed. Above we state (\ref{6.10}) separately just for convenience.
\end{rem}
The third lemma we need is the following large deviation estimate on the extreme eigenvalue of Wigner matrices under the assumption $\mathbf{C_0}$.
\begin{lem} \label{lem.12}
Suppose that $M_n$ is a Wigner matrix obeying the Condition $\mathbf{C_0}$, one has
\begin{eqnarray*}
 \max_{l}|\lambda_l(M_n)|\leq 2+\epsilon
\end{eqnarray*}
 with overwhelming probability for any fixed $\epsilon>0$. Moreover, there exists for any $K\geq 3$
\begin{eqnarray*}
\mathbb{P}(\max_{l}|\lambda_l(M_n)|\geq K)\leq \exp(-cn^c\log K)
\end{eqnarray*}
 for sufficiently large $n$
\end{lem}
\begin{proof} Observe that the first statement is a direct consequence of Lemma {\ref{lem.9}}. And the second one  was proved by Erd\H{o}s, Yau and Yin in \cite{EYY1} (See Lemma 7.2 of \cite{EYY1}).
\end{proof}
The last one is a uniform variance estimate for any single eigenvalue in the bulk, which was proved recently by Dallaporta in \cite{Da}.
 \begin{lem}[eigenvalue variance bounds, \cite{Da}] \label{lem.10}
 Let $M_n$ be a complex Wigner matrices satisfying {\emph{Condition}} $\mathbf{C_0}$, then for any $0<\delta\leq \frac12$, there exists a constant $C(\delta)>0$ such that for all $n\geq 2$ and $\delta n\leq l\leq (1-\delta)n$,
 \begin{eqnarray*}
 Var\{\lambda_l\}\leq C(\delta)\frac{\log n}{n^2}.
 \end{eqnarray*}
 \end{lem}

\end{document}